\documentclass[11pt,twoside]{article}

\usepackage{amsfonts}
\usepackage{amsmath,amssymb,amsthm}
\usepackage{graphicx}
\usepackage{epstopdf}

\usepackage{fancyhdr}
\usepackage{graphics}
\usepackage{amsfonts}
\usepackage{amsmath}

\usepackage{color}

\theoremstyle{plain}

\newtheorem{theo}{Theorem}[section]

\newtheorem{lem}{Lemma}[section]
\newtheorem{prop}{Proposition}[section]
\newtheorem{cor}{Corollary}[section]

\theoremstyle{definition} 

\newtheorem{nota}{Notation}[section]
\newtheorem{de}{Definition}[section]
\newtheorem{exa}{Example}[section]
\newtheorem{as}{Assumption}[section]
\newtheorem{alg}{Algorithm}[section]

\newcommand{\btheo}{\begin{theo}}
\newcommand{\bde}{\begin{de}}
\newcommand{\ble}{\begin{lem}}
\newcommand{\bpr}{\begin{prop}}
\newcommand{\bno}{\begin{nota}}
\newcommand{\bex}{\begin{exa}}
\newcommand{\bcor}{\begin{cor}}
\newcommand{\spro}{\begin{proof}}
\newcommand{\bas}{\begin{as}}
\newcommand{\balg}{\begin{alg}}

\newcommand{\etheo}{\end{theo}}
\newcommand{\ede}{\end{de}}
\newcommand{\ele}{\end{lem}}
\newcommand{\epr}{\end{prop}}
\newcommand{\eno}{\end{nota}}
\newcommand{\eex}{\end{exa}}
\newcommand{\ecor}{\end{cor}}
\newcommand{\fpro}{\end{proof}}
\newcommand{\eas}{\end{as}}
\newcommand{\ealg}{\end{alg}}

\theoremstyle{plain}

\newtheorem{theos}{Theorem}
\newtheorem{props}{Proposition}
\newtheorem{lems}{Lemma}
\newtheorem{cors}{Corollary}

\theoremstyle{definition}
\newtheorem{exas}{Example}
\newtheorem{algs}{Algorithm}
\newtheorem{asss}{Assumption}
\newtheorem{defns}{Definition}

\newcommand{\btheos}{\begin{theos}}
\newcommand{\etheos}{\end{theos}}
\newcommand{\bprops}{\begin{props}}
\newcommand{\eprops}{\end{props}}
\newcommand{\bdes}{\begin{defns}}
\newcommand{\edes}{\end{defns}}
\newcommand{\blems}{\begin{lems}}
\newcommand{\elems}{\end{lems}}
\newcommand{\bcors}{\begin{cors}}
\newcommand{\ecors}{\end{cors}}
\newcommand{\bexs}{\begin{exas}}
\newcommand{\eexs}{\end{exas}}
\newcommand{\balgs}{\begin{algs}}
\newcommand{\ealgs}{\end{algs}}
\newcommand{\bass}{\begin{asss}}
\newcommand{\eass}{\end{asss}}

\newcommand{\rdim}{\ensuremath{r}}
\newcommand{\cov}{\ensuremath{\operatorname{cov}}}

\long\def\comment#1{}

\newcommand{\xls}{\ensuremath{x^{\mbox{\tiny{LS}}}}}
\newcommand{\Xls}{\ensuremath{X^{\mbox{\tiny{LS}}}}}

\newcommand{\numobs}{\ensuremath{n}}
\newcommand{\usedim}{\ensuremath{d}}
\newcommand{\kdim}{\ensuremath{s}}
\newcommand{\plaincon}{\ensuremath{c}}
\newcommand{\xstar}{\ensuremath{x^*}}

\newcommand{\xhat}{\ensuremath{\widehat{x}}}

\newcommand{\mprob}{\ensuremath{\mathbb{P}}}

\newcommand{\numproj}{\ensuremath{m}}

\newcommand{\Sketch}{\ensuremath{S}}
\newcommand{\Width}{\ensuremath{\mathcal{W}}}
\newcommand{\Constraint}{\ensuremath{\mathcal{C}}}
\newcommand{\Amat}{\ensuremath{A}} 
\newcommand{\yvec}{\ensuremath{y}}

\newcommand{\real}{\ensuremath{\mathbb{R}}}
\newcommand{\defn}{\ensuremath{: \, =}}

\newcommand{\inprod}[2]{\ensuremath{\langle #1 , \, #2 \rangle}}

\newcommand{\KCONELS}{\ensuremath{\KCONE^{\mbox{\tiny{LS}}}}}

\newcommand{\Sphere}[1]{\ensuremath{\mathcal{S}^{#1-1}}}

\newcommand{\CEXP}[1]{\ensuremath{e^{#1}}}

\newcommand{\vhat}{\ensuremath{\widehat{v}}}

\newcommand{\Exs}{\ensuremath{\mathbb{E}}}

\newcommand{\diag}{\ensuremath{\operatorname{diag}}}

\newcommand{\fixvec}{\ensuremath{u}}

\newcommand{\PackNum}{\ensuremath{M}}

\newcommand{\ENCMIN}[1]{\ensuremath{\Big \{#1 \Big \}}}

\newcommand{\xit}[1]{\ensuremath{x^{#1}}}
\newcommand{\SketchIt}[1]{\ensuremath{\Sketch^{#1}}}
\newcommand{\ZSUP}{\ensuremath{Z_2}}
\newcommand{\ZINF}{\ensuremath{Z_1}}
\newcommand{\ZSUPIT}[1]{\ensuremath{\ZSUP(\SketchIt{#1})}}
\newcommand{\ZINFIT}[1]{\ensuremath{\ZINF(\SketchIt{#1})}}

\newcommand{\order}{\ensuremath{\mathcal{O}}}

\newlength{\widebarargwidth}
\newlength{\widebarargheight}
\newlength{\widebarargdepth}
\DeclareRobustCommand{\widebar}[1]{%
  \settowidth{\widebarargwidth}{\ensuremath{#1}}%
  \settoheight{\widebarargheight}{\ensuremath{#1}}%
  \settodepth{\widebarargdepth}{\ensuremath{#1}}%
  \addtolength{\widebarargwidth}{-0.3\widebarargheight}%
  \addtolength{\widebarargwidth}{-0.3\widebarargdepth}%
  \makebox[0pt][l]{\hspace{0.3\widebarargheight}%
    \hspace{0.3\widebarargdepth}%
    \addtolength{\widebarargheight}{0.3ex}%
    \rule[\widebarargheight]{0.95\widebarargwidth}{0.1ex}}%
  {#1}}

\newcommand{\Tcrit}{\ensuremath{{N}}}
\newcommand{\Vvar}{\ensuremath{\widebar{Y}}}

\newcommand{\xtil}{\ensuremath{\widetilde{x}}}

\newcommand{\kull}[2]{\ensuremath{D(#1\; \| \; #2)}}
\newcommand{\Prob}{\ensuremath{\mathbb{P}}}

\newcommand{\NORMAL}{\ensuremath{\mathcal{N}}}

\newcommand{\wvec}{\ensuremath{w}}

\newcommand{\SEMI}[1]{\ensuremath{\|#1\|_{\Amat}}}
\newcommand{\SEMIFRO}[1]{\ensuremath{\|#1\|_{\Amat}}}

\newcommand{\KCONE}{\ensuremath{\mathcal{K}}}

\newcommand{\Event}{\ensuremath{\mathcal{E}}}
\newcommand{\ytil}{\ensuremath{\widetilde{\yvec}}}

\newcommand{\Ball}{\ensuremath{\mathbb{B}}}

\newcommand{\goodendex}{\ensuremath{\diamondsuit}}

\newcommand{\uhat}{\ensuremath{\widehat{u}}}

\newcommand{\Ind}{\ensuremath{\mathbb{I}}}

\newcommand{\widgraph}[2]{\includegraphics[keepaspectratio,width=#1]{#2}}

\DeclareMathOperator{\rank}{rank}

\newcommand{\matsnorm}[2]{|\!|\!| #1 | \! | \!|_{{#2}}}

\newcommand{\opnorm}[1]{\ensuremath{\matsnorm{#1}{\mbox{\tiny{op}}}}}
\newcommand{\nucnorm}[1]{\ensuremath{\matsnorm{#1}{\mbox{\tiny{nuc}}}}}
\newcommand{\fronorm}[1]{\ensuremath{\matsnorm{#1}{\mbox{\tiny{fro}}}}}

\newcommand{\NonGaussComp}{\ensuremath{\mathcal{W}}}
\newcommand{\NewPlain}{\ensuremath{\Theta}}
\newcommand{\KCONESTAR}{\ensuremath{\KCONE^*}}
\newcommand{\DelCrit}[1]{\ensuremath{\varepsilon_{#1}}}
\newcommand{\KCONEUN}{\ensuremath{\widebar{\KCONE}}}

\newcommand{\Xstar}{\ensuremath{X^*}}
\newcommand{\Ymat}{\ensuremath{Y}}
\newcommand{\Wmat}{\ensuremath{W}}
\newcommand{\usedima}{\ensuremath{d_1}}
\newcommand{\usedimb}{\ensuremath{d_2}}

\newcommand{\xdagger}{\ensuremath{x^\dagger}}

\newcommand{\EVENTIT}[1]{\ensuremath{\Event^{#1}}}
\newcommand{\deln}{\ensuremath{\varepsilon_\numobs}}
\newcommand{\AuxEvent}{\ensuremath{\mathcal{B}}}

\newcommand{\MyMax}{V_\numobs} 
\newcommand{\starset}{\ensuremath{\operatorname{star}}}
\newcommand{\ztil}{\ensuremath{\widetilde{z}}}
\newcommand{\sketch}{\ensuremath{s}}

\newcommand{\PackNumTwo}[1]{\ensuremath{\PackNum_{#1}}}

\newcommand{\gsketch}{\ensuremath{g_\Sketch}}

\newcommand{\RHODEL}{\ensuremath{\rho}}

\newcommand{\MYEPS}{\ensuremath{\varepsilon}}

\newcommand{\spindex}{\ensuremath{\kdim}}
\newcommand{\COMP}{\ensuremath{C}}

\newcommand{\ntest}{\ensuremath{\numobs_{\mbox{\footnotesize{test}}}}}
\newcommand{\ntrain}{\ensuremath{\numobs_{\mbox{\footnotesize{train}}}}}
\newcommand{\dtask}{\ensuremath{\usedim_{\mbox{\footnotesize{task}}}}}

\newcommand{\tracer}[2]{\ensuremath{\langle \!\langle {#1}, \; {#2}
\rangle \!\rangle}}

\newcommand{\sigval}{\ensuremath{\gamma}}
\newcommand{\gammin}{\ensuremath{\sigval_{\operatorname{\tiny{min}}}}}
\newcommand{\gammax}{\ensuremath{\sigval_{\operatorname{\tiny{max}}}}}

\newcommand{\sign}{\ensuremath{\operatorname{sign}}}

\newcommand{\Atil}{\ensuremath{\widetilde{A}}}
\newcommand{\clconv}{\ensuremath{\operatorname{clconv}}}

\newcommand{\var}{\ensuremath{\operatorname{var}}}
\makeatletter
\long\def\@makecaption#1#2{
        \vskip 0.8ex
        \setbox\@tempboxa\hbox{\small {\bf #1:} #2}
        \parindent 1.5em  
        \dimen0=\hsize
        \advance\dimen0 by -3em
        \ifdim \wd\@tempboxa >\dimen0
                \hbox to \hsize{
                        \parindent 0em
                        \hfil 
                        \parbox{\dimen0}{\def\baselinestretch{0.96}\small
                                {\bf #1.} #2
                                } 
                        \hfil}
        \else \hbox to \hsize{\hfil \box\@tempboxa \hfil}
        \fi
        }
\makeatother

\begin{document}

\begin{center} {\LARGE{\bf{ Iterative Hessian sketch:  Fast and accurate
solution approximation for constrained least-squares}}} \\
  \vspace{1cm}
\begin{tabular}{ccc}
  {\large Mert Pilanci$^{1}$} & & {\large{Martin J.\ Wainwright$^{1,2}$}}
\end{tabular}

 \vspace{.2cm}
  \texttt{\{mert, wainwrig\}@berkeley.edu} \\
  \vspace{1cm}
  {\large University of California, Berkeley} \\
  \vspace{.15cm} $^1$Department of Electrical Engineering and Computer
  Science ~~~~ $^2$Department of Statistics

\vspace*{.2in}

\today

\end{center}

\vspace*{.5in}


\begin{abstract}
We study randomized sketching methods for approximately solving
least-squares problem with a general convex constraint.  The quality
of a least-squares approximation can be assessed in different ways:
either in terms of the value of the quadratic objective function (cost
approximation), or in terms of some distance measure between the
approximate minimizer and the true minimizer (solution approximation).
Focusing on the latter criterion, our first main result provides a
general lower bound on any randomized method that sketches both the
data matrix and vector in a least-squares problem; as a surprising
consequence, the most widely used least-squares sketch is sub-optimal
for solution approximation.  We then present a new method known as the
\emph{iterative Hessian sketch}, and show that it can be used to
obtain approximations to the original least-squares problem using a
projection dimension proportional to the statistical complexity of the
least-squares minimizer, and a logarithmic number of iterations.  We
illustrate our general theory with simulations for both unconstrained
and constrained versions of least-squares, including
$\ell_1$-regularization and nuclear norm constraints. We also
numerically demonstrate the practicality of our approach in a real
face expression classification experiment.
\end{abstract}


\section{Introduction}

Over the past decade, the explosion of data volume and complexity has
led to a surge of interest in fast procedures for approximate forms of
matrix multiplication, low-rank approximation, and convex
optimization.  One interesting class of problems that arise frequently
in data analysis and scientific computing are constrained
least-squares problems.  More specifically, given a data vector $\yvec
\in \real^\numobs$, a data matrix $\Amat \in \real^{\numobs \times
  \usedim}$ and a convex constraint set $\Constraint$, a constrained
least-squares problem can be written as follows
\begin{align}
\label{EqnConstrainedLeastSquares}
\xls & \defn \arg \min_{x \in \Constraint} f(x) \qquad \mbox{where
  $f(x) \defn \frac{1}{2 \numobs} \|\Amat x - \yvec\|_2^2$.}
\end{align}
The simplest case is the unconstrained form ($\Constraint =
\real^\usedim$), but this class also includes other interesting
constrained programs, including those based $\ell_1$-norm balls,
nuclear norm balls, interval constraints $[-1,1]^d$ and other types of regularizers designed to
enforce structure in the solution.

Randomized sketches are a well-established way of obtaining an
approximate solutions to a variety of problems, and there is a long
line of work on their uses (e.g., see the books and
papers~\cite{Vem04, BouDri09,Mahoney11, DriMahMutSar09, KanNel14}, as
well as references therein). In application to
problem~\eqref{EqnConstrainedLeastSquares}, sketching methods
involving using a random matrix \mbox{$\Sketch \in \real^{\numproj
    \times \numobs}$} to project the data matrix $\Amat$ and/or data
vector $\yvec$ to a lower dimensional space ($\numproj \ll \numobs$),
and then solving the approximated least-squares problem.  There are
many choices of random sketching matrices; see
Section~\ref{SecDifferentTypes} for discussion of a few possibilities.
Given some choice of random sketching matrix $\Sketch$, the most
well-studied form of sketched least-squares is based on solving the
problem
\begin{align}
\label{EqnSuboptimalSketch}
\xtil & \defn \arg \min_{x \in \Constraint} \ENCMIN{\frac{1}{2
    \numobs} \| \Sketch \Amat x - \Sketch \yvec\|_2^2},
\end{align}
in which the data matrix-vector pair $(\Amat, \yvec)$ are approximated
by their sketched versions $(\Sketch \Amat, \Sketch \yvec)$.  Note
that the sketched program is an $\numproj$-dimensional least-squares
problem, involving the new data matrix $\Sketch \Amat \in
\real^{\numproj \times \usedim}$.  Thus, in the regime $\numobs \gg
\usedim$, this approach can lead to substantial computational savings
as long as the projection dimension $\numproj$ can be chosen
substantially less than $\numobs$.  A number of authors
(e.g.,~\cite{BouDri09,DriMahMutSar09,Mahoney11,PilWai14a}) have
investigated the properties of this sketched
solution~\eqref{EqnSuboptimalSketch}, and accordingly, we refer to to
it as the \emph{classical least-squares sketch}.

There are various ways in which the quality of the approximate
solution $\xtil$ can be assessed.  One standard way is in terms of the
minimizing value of the quadratic cost function $f$ defining the
original problem~\eqref{EqnConstrainedLeastSquares}, which we refer to
as \emph{cost approximation}. In terms of $f$-cost, the approximate
solution $\xtil$ is said to be $\MYEPS$-optimal if
\begin{align}
\label{EqnCostDeltaOptimal}
f(\xls) \; \leq \; f(\xtil) \; \leq \; (1+\MYEPS)^2 f(\xls).
\end{align}
For example, in the case of unconstrained least-squares ($\Constraint
= \real^\usedim$) with $\numobs > \usedim$, it is known that with
Gaussian random sketches, a sketch size $\numproj \succsim
\frac{1}{\MYEPS^2} \usedim$ suffices to guarantee that $\xtil$ is
\mbox{$\MYEPS$-optimal} with high probability (for instance, see the
papers by Sarlos~\cite{Sarlos2006} and Mahoney~\cite{Mahoney11}, as
well as references therein).  Similar guarantees can be established
for sketches based on sampling according to the statistical leverage
scores~\cite{Drineas2010,Drineas2012fast}.  Sketching can also be
applied to problems with constraints: Boutsidis and
Drineas~\cite{BouDri09} prove analogous results for the case of
non-negative least-squares considering the sketch in \eqref{EqnSuboptimalSketch}, whereas our own past work~\cite{PilWai14a}
provides sufficient conditions for $\MYEPS$-accurate cost
approximation of least-squares problems over arbitrary convex sets based also on the form in \eqref{EqnSuboptimalSketch}.

It should be noted, however, that other notions of ``approximation
goodness'' are possible.  In many applications, it is the
least-squares minimizer $\xls$ itself---as opposed to the cost value
$f(\xls)$---that is of primary interest.  In such settings, a more
suitable measure of approximation quality would be the $\ell_2$-norm
$\|\xtil - \xls\|_2$, or the prediction (semi)-norm
\begin{align}
\label{EqnPredSemiNorm}
\SEMI{ \xtil - \xls} \defn \frac{1}{\sqrt{\numobs}} \|\Amat(\xtil -
\xls)\|_2.
\end{align}
We refer to these measures as \emph{solution approximation}.

Now of course, a cost approximation bound~\eqref{EqnCostDeltaOptimal}
can be used to derive guarantees on the solution approximation error.
However, it is natural to wonder whether or not, for a reasonable
sketch size, the resulting guarantees are ``good''. For instance,
using arguments from Drineas et al.~\cite{DriMahMutSar09}, for the
problem of unconstrained least-squares, it can be shown that the
same conditions ensuring a $\MYEPS$-accurate cost approximation also
ensure that
\begin{align}
\label{EqnClassicalBound}
\SEMI{\xtil - \xls} & \leq \MYEPS \, \sqrt{f(\xls)}.
\end{align}
Given lower bounds on the singular values of the data matrix $\Amat$,
this bound also yields control of the $\ell_2$-error.

In certain ways, the bound~\eqref{EqnClassicalBound} is quite
satisfactory: given our normalized
definition~\eqref{EqnConstrainedLeastSquares} of the least-squares
cost $f$, the quantity $f(\xls)$ remains an order one quantity as the
sample size $\numobs$ grows, and the multiplicative factor $\MYEPS$
can be reduced by increasing the sketch dimension $\numproj$.  But how
small should $\MYEPS$ be chosen?  In many applications of
least-squares, each element of the response vector $y \in
\real^\numobs$ corresponds to an observation, and so as the sample
size $\numobs$ increases, we expect that $\xls$ provides a more
accurate approximation to some underlying population quantity, say
$\xstar \in \real^\usedim$.  As an illustrative example, in the
special case of unconstrained least-squares, the accuracy of the
least-squares solution $\xls$ as an estimate of $\xstar$ scales as
$\SEMI{\xls - \xstar} \asymp \frac{\sigma^2 \usedim}{\numobs}$.
Consequently, in order for our sketched solution to have an accuracy
of the same order as the least-square estimate, we must set $\MYEPS^2
\asymp \frac{\sigma^2 \usedim}{\numobs}$.  Combined with our earlier
bound on the projection dimension, this calculation suggests that a
projection dimension of the order
\begin{align*}
\numproj & \succsim \frac{\usedim}{\MYEPS^2} \; \asymp
\frac{\numobs}{\sigma^2}
\end{align*}
is required.  This scaling is undesirable in the regime $\numobs \gg
\usedim$, where the whole point of sketching is to have the sketch
dimension $\numproj$ much lower than $\numobs$.\\

\begin{figure}[h]
\begin{center}
\begin{tabular}{cc}
\widgraph{.45 \textwidth}{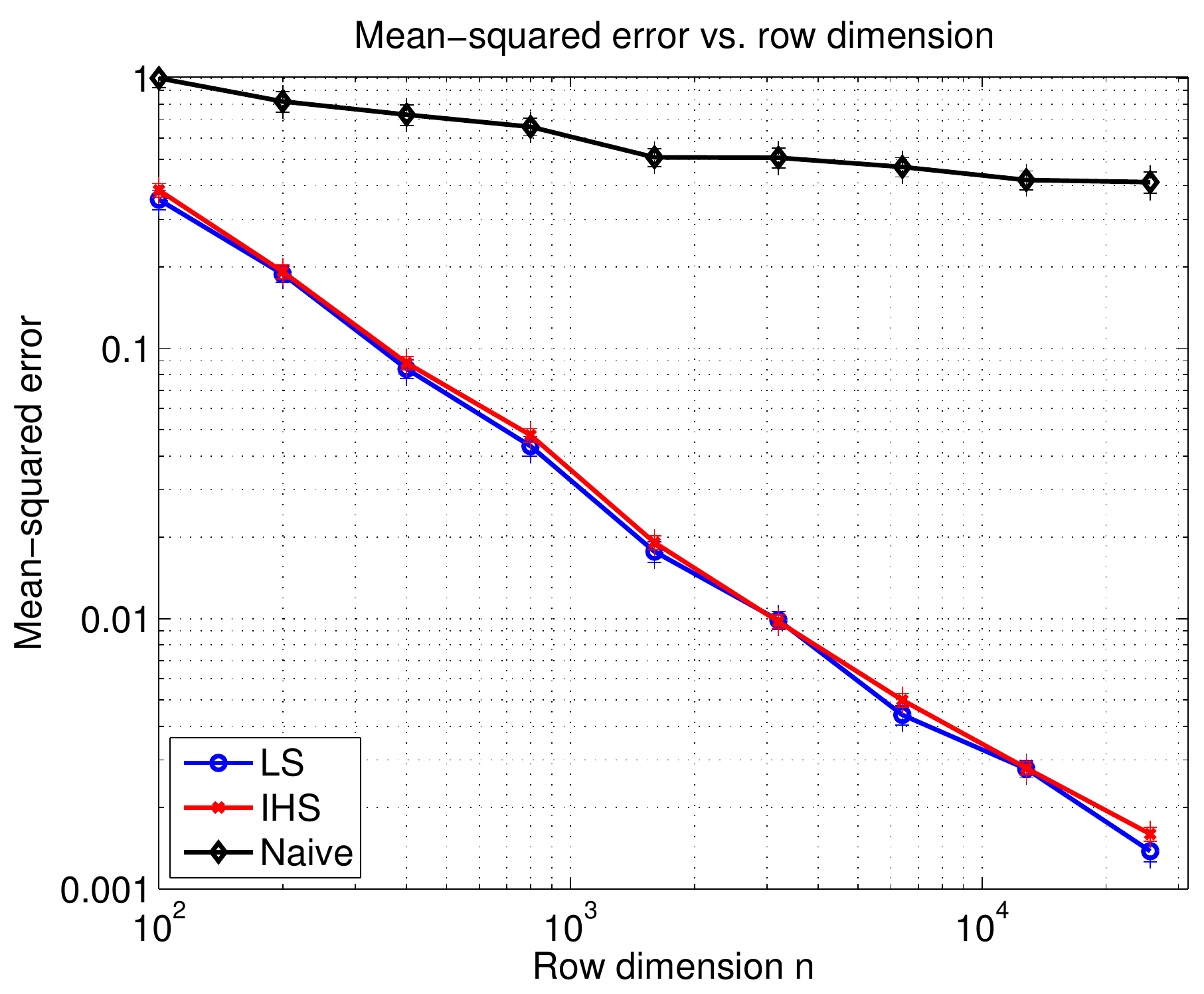} & 
\widgraph{.45 \textwidth}{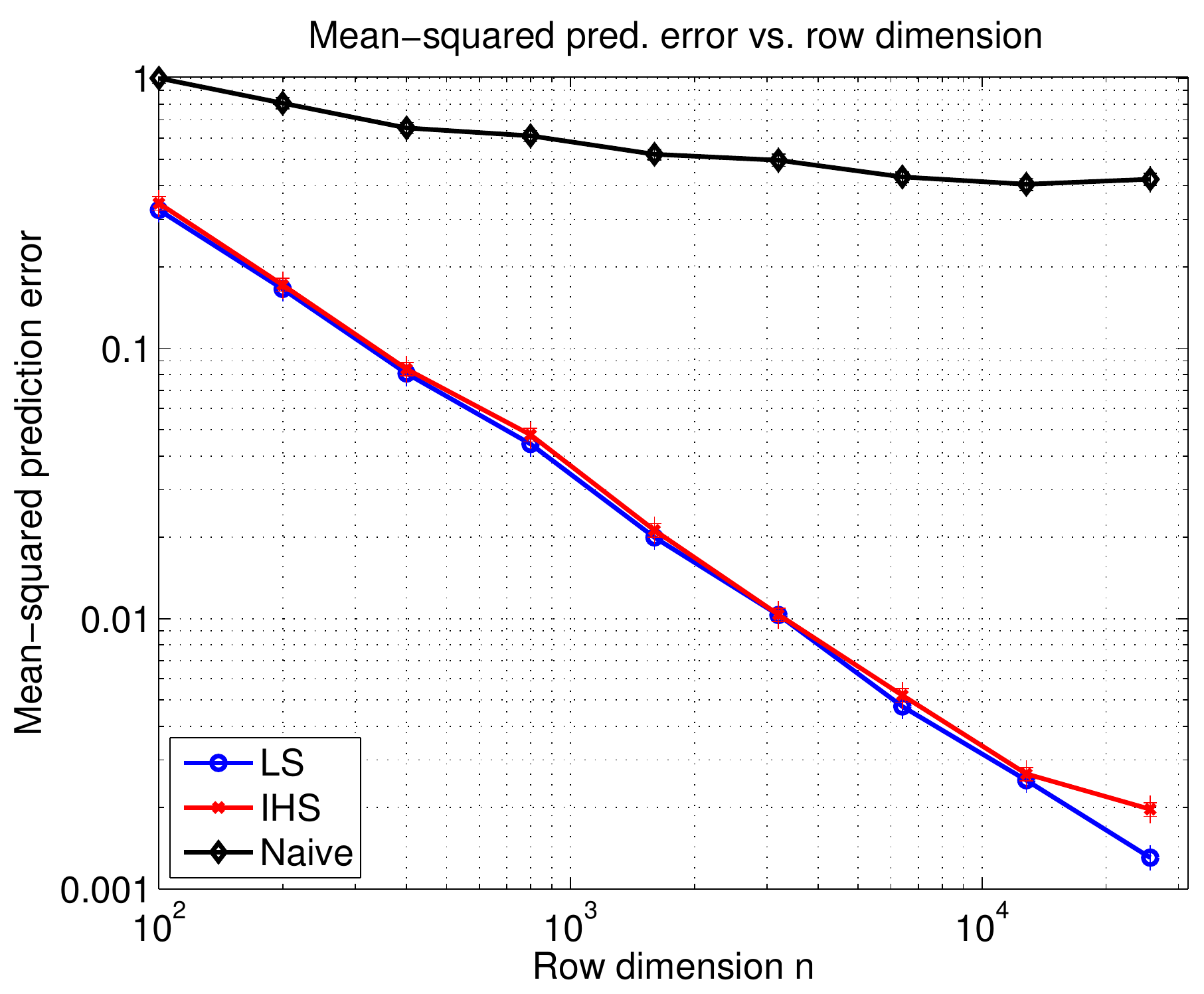} \\
(a) & (b)
\end{tabular}
\caption{Plots of mean-squared error versus the row dimension $\numobs
  \in \{100, 200, 400, \ldots, 25600 \}$ for unconstrained
  least-squares in dimension $\usedim = 10$.  The blue curves
  correspond to the error $\xls - \xstar$ of the unsketched
  least-squares estimate. Red curves correspond to the IHS method
  applied for $\Tcrit = 1 + \lceil \log(\numobs) \rceil$ rounds using
  a sketch size $\numproj = 7 \usedim$.  Black curves correspond to
  the naive sketch applied using $M = \Tcrit \numproj$ projections in
  total, corresponding to the same number used in all iterations of
  the IHS algorithm. (a) Error $\|\xtil - \xstar\|_2^2$.  (b)
  Prediction error $\SEMI{\xtil - \xstar}^2 = \frac{1}{\numobs}
  \|\Amat (\xtil - \xstar)\|_2^2$.  Each point corresponds to the mean
  taken over $300$ trials with standard errors shown above and below
  in crosses.}
\label{FigLeastSquaresFixedM}
\end{center}
\end{figure}

Now the alert reader will have observed that the preceding argument
was only rough and heuristic.  However, the first result of this paper
(Theorem~\ref{ThmSubOptimal}) provides a rigorous confirmation of the
conclusion: whenever $\numproj \ll \numobs$, the classical
least-squares sketch~\eqref{EqnSuboptimalSketch} is sub-optimal as a
method for solution approximation. Figure~\ref{FigLeastSquaresFixedM}
provides an empirical demonstration of the poor behavior of the
classical least-squares sketch for an unconstrained problem.  

This sub-optimality holds not only for unconstrained least-squares but
also more generally for a broad class of constrained problems.
Actually, Theorem~\ref{ThmSubOptimal} is a more general claim:
\emph{any estimator} based only on the pair $(\Sketch \Amat, \Sketch
\yvec)$---an infinite family of methods including the standard
sketching algorithm as a particular case---is sub-optimal relative to
the original least-squares estimator in the regime $\numproj \ll
\numobs$.  We are thus led to a natural question: can this
sub-optimality be avoided by a different type of sketch that is
nonetheless computationally efficient?  Motivated by this question,
our second main result (Theorem~\ref{ThmOptIterative}) is to propose
an alternative method---known as the iterative Hessian sketch---and
prove that it yields optimal approximations to the least-squares
solution using a projection size that scales with the intrinsic
dimension of the underlying problem, along with a logarithmic number of iterations. The main idea underlying iterative Hessian sketch is to obtain multiple sketches of the data $(S^1 A, ..., S^{N} A)$ and iteratively refine the solution where $N$ can be chosen logarithmic in $n$.

The remainder of this paper is organized as follows.  In
Section~\ref{SecMain}, we begin by introducing some background on
classes of random sketching matrices, before turning to the statement
of our lower bound (Theorem~\ref{ThmSubOptimal}) on the classical
least-squares sketch~\eqref{EqnSuboptimalSketch}.  We then introduce
the Hessian sketch, and show that an iterative version of it can be
used to compute $\varepsilon$-accurate solution approximations using
$\log(1/\varepsilon)$-steps (Theorem~\ref{ThmOptIterative}).  In
Section~\ref{SecConsequences}, we illustrate the consequences of this
general theorem for various specific classes of least-squares
problems, and we conclude with a discussion in
Section~\ref{SecDiscussion}.  The majority of our proofs are deferred
to the appendices.

\paragraph{Notation:}  For the convenience of the reader, we summarize
some standard notation used in this paper.  For sequences
$\{a_t\}_{t=0}^\infty$ and $\{b_t\}_{t=0}^\infty$, we use the notation
$a_t \preceq b_t$ to mean that there is a constant (independent of
$t$) such that $a_t \leq C \, b_t$ for all $t$.  Equivalently, we
write $b_t \succeq a_t$.  We write $a_t \asymp b_t$ if $a_t \preceq
b_t$ and $b_t \preceq a_t$.

\section{Main results}
\label{SecMain}

In this section, we begin with background on different classes of
randomized sketches, including those based on random matrices with
sub-Gaussian entries, as well as those based on randomized orthonormal
systems and random sampling.  In Section~\ref{SecSubOpt}, we prove a
general lower bound on the solution approximation accuracy of any
method that attempts to approximate the least-squares problem based on
observing only the pair $(\Sketch \Amat, \Sketch \yvec)$.  This
negative result motivates the investigation of alternative sketching
methods, and we begin this investigation by introducing the Hessian
sketch in Section~\ref{SecHessianSketch}.  It serves as the basic
building block of the iterative Hessian sketch (IHS), which can be
used to construct an iterative method that is optimal up to
logarithmic factors.


\subsection{Different types of randomized sketches}
\label{SecDifferentTypes}

Various types of randomized sketches are possible, and we describe a
few of them here.  Given a sketching matrix $\Sketch$, we use
$\{\sketch_i\}_{i=1}^\numproj$ to denote the collection of its
$\numobs$-dimensional rows.  We restrict our attention to sketch
matrices that are zero-mean, and that are normalized so that
$\Exs[\Sketch^T \Sketch/\numproj] = I_\numobs$.

\paragraph{Sub-Gaussian sketches:}  The most classical sketch
is based on a random matrix \mbox{$\Sketch \in \real^{\numproj \times
    \numobs}$} with i.i.d. standard Gaussian entries.  A
straightforward generalization is a random sketch with
i.i.d. sub-Gaussian rows.  In particular, a zero-mean random vector
$\sketch \in \real^{\numobs}$ is $1$-sub-Gaussian if for any $u \in
\real^\numobs$, we have
\begin{align}
\mprob[ \inprod{\sketch}{u} \geq \MYEPS \|u\|_2 \big] & \leq
\CEXP{-\MYEPS^2/2} \qquad \mbox{for all $\MYEPS \geq 0$.}
\end{align}
For instance, a vector with i.i.d. $N(0,1)$ entries is
$1$-sub-Gaussian, as is a vector with i.i.d. Rademacher entries
(uniformly distributed over $\{-1, +1\}$).  Suppose that we generate a
random matrix $\Sketch \in \real^{\numproj \times \numobs}$ with
i.i.d.  rows that are zero-mean, $1$-sub-Gaussian, and with
$\cov(\sketch) = I_\numobs$; we refer to any such matrix as a
\emph{sub-Gaussian sketch}.  As will be clear, such sketches are the
most straightforward to control from the probabilistic point of view.
However, from a computational perspective, a disadvantage of
sub-Gaussian sketches is that they require matrix-vector
multiplications with unstructured random matrices.  In particular,
given an data matrix $\Amat \in \real^{\numobs \times \usedim}$,
computing its sketched version $\Sketch \Amat$ requires
$\order(\numproj \numobs \usedim)$ basic operations in general (using classical matrix multiplication).


\paragraph{Sketches based on randomized orthonormal systems (ROS):}

The second type of randomized sketch we consider is \emph{randomized
  orthonormal system} (ROS), for which matrix multiplication can be
performed much more efficiently.

In order to define a ROS sketch, we first let $H \in \real^{\numobs
  \times \numobs}$ be an orthonormal matrix with entries $H_{ij} \in [
  -\frac{1}{\sqrt{\numobs}}, \frac{1}{\sqrt{\numobs}} ]$.  Standard
classes of such matrices are the Hadamard or Fourier bases, for which
matrix-vector multiplication can be performed in $\order(\numobs \log
\numobs)$ time via the fast Hadamard or Fourier transforms,
respectively.  Based on any such matrix, a sketching matrix $\Sketch
\in \real^{\numproj \times \numobs}$ from a ROS ensemble is obtained
by sampling i.i.d. rows of the form
\begin{align*}
\sketch^T & = \sqrt{\numobs} e_j^T H D \qquad \mbox{with probability
  $1/\numobs$ for $j = 1, \ldots, \numobs$},
\end{align*}
where the random vector $e_j \in \real^\numobs$ is chosen uniformly at
random from the set of all $\numobs$ canonical basis vectors, and $D =
\diag(\nu)$ is a diagonal matrix of i.i.d. Rademacher variables
\mbox{$\nu \in \{-1, +1\}^\numobs$.}  Given a fast routine for
matrix-vector multiplication, the sketched data $(\Sketch A, \Sketch
y)$ can be formed in $\order(\numobs \, \usedim \log \numproj)$ time
(for instance, see the paper~\cite{ailon2006approximate}).

\paragraph{Sketches based on random row sampling:}
Given a probability distribution $\{p_j\}_{j=1}^\numobs$ over
$[\numobs] = \{1, \ldots, \numobs \}$, another choice of sketch is to
randomly sample the rows of the extended data matrix $\begin{bmatrix}
  \Amat & \yvec
\end{bmatrix}$ a total of $\numproj$ times with replacement
from the given probability distribution.  Thus, the rows of $S$ are
independent and take on the values
\begin{align*}
s^T & = \frac{e_j}{\sqrt{p_j}} \qquad \mbox{with probability $p_j$ for
  $j = 1, \ldots, \numobs$}
\end{align*}
where $e_j \in \real^\numobs$ is the $j^{th}$ canonical basis vector.
Different choices of the weights $\{p_j\}_{j=1}^\numobs$ are possible,
including those based on the leverage values of $\Amat$---i.e., $p_j
\propto \|u_j\|_2$ for $j= 1, \ldots, \numobs$, where $U \in
\real^{\numobs \times \usedim}$ is the matrix of left singular vectors
of $\Amat$ \cite{Drineas2010}.  In our analysis of lower bounds to follow, we assume that
the weights are $\alpha$-balanced, meaning that
\begin{align}
\label{EqnAlphaBalanced}
\max_{j=1, \ldots, \numobs} p_j & \leq \frac{\alpha}{\numobs}
\end{align}
for some constant $\alpha$ independent of $\numobs$.

In the following section, we present a lower bound that applies to all
the three kinds of sketching matrices described above.


\subsection{Sub-optimality of classical least-squares sketch}
\label{SecSubOpt}

We begin by proving a lower bound on any estimator that is a function
of the pair $(\Sketch \Amat, \Sketch \yvec)$.  In order to do so,
we consider an ensemble of least-squares problems, namely those
generated by a noisy observation model of the form 
\begin{align}
\yvec & = \Amat \xstar + \wvec, \qquad \mbox{where $\wvec \sim N(0,
  \sigma^2 I_{\numobs})$,}
\end{align}
the data matrix $\Amat \in \real^{\numobs \times \usedim}$ is fixed,
and the unknown vector $\xstar$ belongs to some compact subset
$\Constraint_0 \subseteq \Constraint$.  In this case, the constrained
least-squares estimate $\xls$ from
equation~\eqref{EqnConstrainedLeastSquares} corresponds to a
constrained form of maximum-likelihood for estimating the unknown
regression vector $\xstar$.  In Appendix~\ref{AppPropMLE}, we provide
a general upper bound on the error $\Exs[ \SEMI{\xls - \xstar}^2]$ in
the least-squares solution as an estimate of $\xstar$.  This result
provides a baseline against which to measure the performance of a
sketching method: in particular, our goal is to characterize the
minimal projection dimension $\numproj$ required in order to return an
estimate $\xtil$ with an error guarantee \mbox{$\SEMI{\xtil - \xls}
  \approx \SEMI{\xls - \xstar}$.}  The result to follow shows that
unless $\numproj \geq \numobs$, then \emph{any method} based on
observing \emph{only} the pair $(\Sketch \Amat, \Sketch \yvec)$
necessarily has a substantially larger error than the least-squares
estimate.  In particular, our result applies to an arbitrary
measureable function $(\Sketch \Amat, \Sketch \yvec) \mapsto
\xdagger$, which we refer to as an \emph{estimator}.

More precisely, our lower bound applies to any random matrix $\Sketch
\in \real^{\numproj \times \numobs}$ for which
\begin{align}
\label{EqnHanaShouldSleep}
\opnorm{\Exs \Big[ \Sketch^T (\Sketch \Sketch^T)^{-1}\Sketch \Big]}
\leq \eta \; \frac{\numproj}{\numobs},
 \end{align}
where $\eta$ is a constant independent of $\numobs$ and $\numproj$,
and $\opnorm{A}$ denotes the $\ell_2$-operator norm (maximum
eigenvalue for a symmetric matrix).  In
Appendix~\ref{AppLemHanaShouldSleep}, we show that these conditions
hold for various standard choices, including most of those discussed
in the previous section.  Our lower bound also involves the complexity
of the set $\Constraint_0$, which we measure in terms of its metric
entropy. In particular, for a given semi-norm $\|\cdot\|$ and
tolerance $\delta > 0$, the $\delta$-packing number
$\PackNumTwo{\delta}$ of the set $\Constraint_0$ is the largest number
of vectors $\{x^j\}_{j=1}^M \subset \Constraint_0$ such that $\|x^j -
x^k \| > \delta$ for all distinct pairs $j \neq k$.

With this set-up, we have the following result: \\

\btheos[Sub-optimality]
\label{ThmSubOptimal}
For any random sketching matrix $\Sketch \in \real^{\numproj \times
  \numobs}$ satisfying condition~\eqref{EqnHanaShouldSleep}, any
estimator $(\Sketch \Amat, \Sketch \yvec) \mapsto \xdagger$ has MSE
lower bounded as
\begin{align}
\label{EqnMinimax}
\sup_{\xstar \in \Constraint_0} \Exs_{\Sketch, \wvec} \big[
  \SEMI{\xdagger - \xstar}^2 \big] & \geq \frac{\sigma^2}{128 \, \eta}
\; \frac{\log(\frac{1}{2}M_{1/2})}{\min\{\numproj, \numobs\}}
\end{align}
where $M_{1/2}$ is the $1/2$-packing number of $\Constraint_0$ in the
semi-norm $\SEMI{\cdot}$.
\etheos
\noindent The proof, given in Appendix~\ref{AppThmSubOptimal}, is
based on a reduction from statistical minimax theory combined with
information-theoretic bounds.  The lower bound is best understood by
considering some concrete examples: \\

\bexs[Sub-optimality for ordinary least-squares]
\label{ExaOrdinaryLeastSquares}
We begin with the simplest case---namely, in which \mbox{$\Constraint
  = \real^\usedim$.}  With this choice and for any data matrix $\Amat$
with $\rank(\Amat) = \usedim$, it is straightforward to show that the
least-squares solution $\xls$ has its prediction mean-squared error at
most
\begin{subequations}
\begin{align}
\label{EqnUpperOLS}
\Exs \big[ \SEMI{\xls - \xstar}^2 \big] & \precsim \: \frac{\sigma^2
  \usedim}{\numobs}.
\end{align}
On the other hand, with the choice $\Constraint_0 = \Ball_2(1)$, we
can construct a $1/2$-packing with $\PackNum = 2^{\usedim}$ elements,
so that Theorem~\ref{ThmSubOptimal} implies that any estimator
$\xdagger$ based on $(\Sketch \Amat, \Sketch \yvec)$ has its
prediction MSE lower bounded as
\begin{align}
\label{EqnLowerOLS}
\Exs_{\Sketch, \wvec} \big[ \SEMI{\xhat - \xstar}^2 \big] & \succsim
\frac{\sigma^2 \, \usedim}{\min\{\numproj, \numobs\}}.
\end{align}
\end{subequations}
Consequently, the sketch dimension $\numproj$ must grow proportionally
to $\numobs$ in order for the sketched solution to have a mean-squared
error comparable to the original least-squares estimate.  This is
highly undesirable for least-squares problems in which $\numobs \gg
\usedim$, since it should be possible to sketch down to a dimension
proportional to $\rank(\Amat) = \usedim$.  Thus,
Theorem~\ref{ThmSubOptimal} this reveals a surprising gap between the
classical least-squares sketch~\eqref{EqnSuboptimalSketch} and the
accuracy of the original least-squares estimate.

In contrast, the sketching method of this paper, known as iterative
Hessian sketching (IHS), matches the optimal mean-squared error using
a sketch of size $\usedim + \log(\numobs)$ in each round, and a total
of $\log(\numobs)$ rounds; see Corollary~\ref{CorOptLS} for a precise
statement.  The red curves in Figure~\ref{FigLeastSquaresFixedM} show
that the mean-squared errors ($\|\xhat - \xstar\|_2^2$ in panel (a),
and $\SEMI{\xhat - \xstar}^2$ in panel (b)) of the IHS method using
this sketch dimension closely track the associated errors of the full
least-squares solution (blue curves).  Consistent with our previous
discussion, both curves drop off at the $\numobs^{-1}$ rate.  

Since the IHS method with $\log(\numobs)$ rounds uses a total of $T =
\log(\numobs) \big \{ \usedim + \log(\numobs) \}$ sketches, a fair
comparison is to implement the classical method with $T$ sketches in
total.  The black curves show the MSE of the resulting sketch: as
predicted by our theory, these curves are relatively flat as a
function of sample size $\numobs$.  Indeed, in this particular case,
the lower bound~\eqref{EqnMinimax}
\begin{align*}
\Exs_{\Sketch, \wvec} \big[ \SEMI{\xtil - \xstar}^2 \big] & \succsim
\frac{\sigma^2 \usedim}{\numproj} \succsim
\frac{\sigma^2}{\log^2(\numobs)},
\end{align*}
showing we can expect (at best) an inverse logarithmic drop-off.
\hfill \goodendex \eexs

\noindent This sub-optimality can be extended to other forms of
constrained least-squares estimates as well, such as those involving
sparsity constraints.
\bexs[Sub-optimality for sparse linear models]
\label{ExaSparseLinear}
We now consider the sparse variant of the linear regression problem,
which involves the $\ell_0$-``ball'' 
\begin{align*}
\Ball_0(\kdim) \defn \big \{ x \in \real^\usedim \mid
\sum_{j=1}^\usedim\Ind[x_j \neq 0] \leq \kdim \},
\end{align*}
corresponding to the set of all vectors with at most $\kdim$ non-zero
entries.  Fixing some radius $R \geq \sqrt{\kdim}$, consider a vector
\mbox{$\xstar \in \Constraint_0 \defn \Ball_0(\kdim) \cap \{\|x\|_1 =
  R \}$,} and suppose that we make noisy observations of the form
$\yvec = \Amat \xstar + \wvec$.

Given this set-up, one way in which to estimate $\xstar$ is by by
computing the least-squares estimate $\xls$ constrained\footnote{This
  set-up is slightly unrealistic, since the estimator is assumed to
  know the radius $R = \|\xstar\|_1$.  In practice, one solves the
  least-squares problem with a Lagrangian constraint, but the
  underlying arguments are basically the same.}  to the $\ell_1$-ball
$\Constraint = \{ x \in \real^\numobs \, \mid \, \|x\|_1 \leq R \}$.
This estimator is a form of the Lasso~\cite{Tibshirani96}: as shown in
Appendix~\ref{AppExamples}, when the design matrix $\Amat$ satisfies
the restricted isometry property (see \cite{CandesTao05} for a definition), then it has MSE at most
\begin{subequations}
\begin{align}
\label{EqnUpperSparse}
\Exs \big[ \SEMI{\xls - \xstar}^2 \big] & \precsim \frac{\sigma^2
  \kdim \log \big( \frac{e \usedim}{\kdim} \big)}{\numobs}.
\end{align}

On the other hand, the $\frac{1}{2}$-packing number $\PackNum$ of the
set $\Constraint_0$ can be lower bounded as \mbox{$\log \PackNum
  \succsim \kdim \log \big(\frac{e \usedim}{\kdim} \big)$;} see
Appendix~\ref{AppExamples} for the details of this calculation.
Consequently, in application to this particular problem,
Theorem~\ref{ThmSubOptimal} implies that any estimator $\xdagger$
based on the pair $(\Sketch \Amat, \Sketch \yvec)$ has mean-squared
error lower bounded as
\begin{align}
\label{EqnLowerSparse}
\Exs_{\wvec, \Sketch} \big[ \SEMI{\xdagger - \xstar}^2 \big] &
\succsim \frac{ \sigma^2 \kdim \log \big(\frac{e \usedim}{\kdim}
  \big)}{\min \{ \numproj, \numobs \}}.
\end{align}
\end{subequations}
Again, we see that the projection dimension $\numproj$ must be of the
order of $\numobs$ in order to match the mean-squared error of the
constrained least-squares estimate $\xls$ up to constant factors.  By
contrast, in this special case, the sketching method developed in this
paper matches the error $\|\xls - \xstar\|_2$ using a sketch dimension
that scales only as $\kdim \log \big( \frac{e \usedim}{\kdim} \big) +
\log(\numobs)$; see Corollary~\ref{CorOptSparse} for the details of a
more general result.
\hfill \goodendex
\eexs

\bexs[Sub-optimality for low-rank matrix estimation]
\label{ExaMultiTask}
In the problem of multivariate regression, the goal is to estimate a
matrix $\Xstar \in \real^{\usedima \times \usedimb}$ model based on
observations of the form
\begin{align}
\Ymat & = \Amat \Xstar + \Wmat,
\end{align}
where $\Ymat \in \real^{\numobs \times \usedima}$ is a matrix of
observed responses, $\Amat \in \real^{\numobs \times \usedima}$ is a
data matrix, and \mbox{$\Wmat \in \real^{\numobs \times \usedimb}$} is
a matrix of noise variables.  One interpretation of this model is as a
collection of $\usedimb$ regression problems, each involving a
$\usedima$-dimensional regression vector, namely a particular column
of $\Xstar$.  In many applications, among them reduced rank
regression, multi-task learning and recommender systems
(e.g.,~\cite{SreAloJaa05,YuaLi06, NegWai09,BunSheWeg11}), it is
reasonable to model the matrix $\Xstar$ as having a low-rank.  Note a
rank constraint on matrix $X$ be written as an $\ell_0$-``norm''
constraint on its singular values: in particular, we have
\begin{align*}
\rank(X) \leq \rdim \quad \mbox{if and only if} \quad \sum_{j=1}^{\min
  \{\usedima, \usedimb \}} \Ind[\sigval_j(X) > 0] \leq \rdim,
\end{align*}
where $\sigval_j(X)$ denotes the $j^{th}$ singular value of $X$.  This
observation motivates a standard relaxation of the rank constraint
using the nuclear norm $\nucnorm{X} \defn \sum_{j=1}^{\min \{\usedima,
  \usedimb \}} \sigval_j(X)$.

Accordingly, let us consider the constrained least-squares problem
\begin{align}
\Xls & = \arg \min_{X \in \real^{\usedima \times \usedimb}}
\ENCMIN{\frac{1}{2} \fronorm{Y - AX}^2} \qquad \mbox{such that
  $\nucnorm{X} \leq R$,}
\end{align}
where $\fronorm{\cdot}$ denotes the Frobenius norm on matrices, or
equivalently the Euclidean norm on its vectorized version.  Let
$\Constraint_0$ denote the set of matrices with rank $\rdim <
\frac{1}{2} \min \{ \usedima, \usedimb \}$, and Frobenius norm at most
one.  In this case, we show in Appendix~\ref{AppPropMLE} that the
constrained least-squares solution $\Xls$ satisfies the bound
\begin{subequations}
\begin{align}
\label{EqnUpperLowRank}
\Exs \Big[ \SEMI{\Xls - \Xstar}^2 \Big] & \precsim \frac{\sigma^2
  \rdim \, (\usedima + \usedimb)}{\numobs}.
\end{align}
On the other hand, the $\frac{1}{2}$-packing number of the set
$\Constraint_0$ is lower bounded as \mbox{$\log \PackNum \succsim
  \rdim \big(\usedima + \usedimb \big)$,} so that
Theorem~\ref{ThmSubOptimal} implies that any estimator $X^\dagger$
based on the pair $(\Sketch \Amat, \Sketch Y)$ has MSE lower bounded
as
\begin{align}
\label{EqnLowerLowRank}
\Exs_{\wvec, \Sketch} \big[ \SEMI{X^\dagger - \Xstar}^2 \big] &
\succsim \frac{ \sigma^2 \rdim \big( \usedima + \usedimb \big)}{\min
  \{ \numproj, \numobs \}}.
\end{align}
\end{subequations}
As with the previous examples, we see the sub-optimality of the
sketched approach in the regime $\numproj < \numobs$.  In contrast,
for this class of problems, our sketching method matches the error
$\SEMIFRO{\Xls - \Xstar}$ using a sketch dimension that scales only as
$\{ \rdim (\usedima + \usedimb) + \log(\numobs) \} \, \log(\numobs$).
See Corollary~\ref{CorLowRank} for further details.

\hfill \goodendex
\eexs


\subsection{Introducing the Hessian sketch}
\label{SecHessianSketch}

As will be revealed during the proof of Theorem~\ref{ThmSubOptimal},
the sub-optimality is in part due to sketching the response
vector---i.e., observing $\Sketch \yvec$ instead of $\yvec$.  It is
thus natural to consider instead methods that sketch \emph{only} the
data matrix $\Amat$, as opposed to both the data matrix and data
vector $y$.  In abstract terms, such methods are based on observing
the pair \mbox{$\big( \Sketch \Amat, \Amat^T \yvec \big) \in
  \real^{\numproj \times \usedim} \times \real^{\usedim}$.} One such
approach is what we refer to as the \emph{Hessian sketch}---namely,
the sketched least-squares problem
\begin{align}
\label{EqnHessianSketch}
\xhat & \defn \arg \min_{x \in \Constraint } \ENCMIN{
  \underbrace{\frac{1}{2} \|\Sketch \Amat x\|_2^2 - \inprod{\Amat^T
      \yvec}{x}}_{\gsketch(x)}}.
\end{align}
As with the classical least-squares
sketch~\eqref{EqnSuboptimalSketch}, the quadratic form is defined by
the matrix $\Sketch \Amat \in \real^{\numproj \times \usedim}$, which
leads to computational savings.  Although the Hessian sketch on its
own does not provide an optimal approximation to the least-squares
solution, it serves as the building block for an iterative method that
can obtain an $\varepsilon$-accurate solution approximation in
$\log(1/\varepsilon)$ iterations.

In controlling the error with respect to the least-squares solution
$\xls$ the set of possible descent directions $\{ x - \xls \, \mid \,
x \in \Constraint \}$ plays an important role.  In particular, we
define the \emph{transformed tangent cone}
\begin{align}
\label{EqnTransformedCone}
\KCONELS & = \big \{ v \in \real^\usedim \, \mid v = t \, \Amat (x -
\xls) \quad \mbox{for some $t \geq 0$ and $x \in \Constraint$} \big
\}.
\end{align}
Note that the error vector $\vhat \defn \Amat (\xhat - \xls)$ of
interest belongs to this cone.  Our approximation bound is a function
of the quantities
\begin{subequations}
\begin{align}
\label{EqnDefnZinf}
\ZINF(\Sketch) & \defn \inf_{v \in \KCONELS \cap \Sphere{\numobs}}
\frac{1}{\numproj} \| \Sketch v\|_2^2 \quad \mbox{and} \\
\label{EqnDefnZsup}
\ZSUP(\Sketch) & \defn \sup_{v \in \KCONELS \cap \Sphere{\numobs}}
\Big| \inprod{\fixvec}{(\frac{\Sketch^T \Sketch}{\numproj} -
  I_\numobs) \, v} \Big|,
\end{align}
\end{subequations}
where $u$ is a fixed unit-norm vector.  These variables played an
important role in our previous analysis~\cite{PilWai14a} of the
classical sketch~\eqref{EqnSuboptimalSketch}. The following bound
applies in a deterministic fashion to any sketching matrix.
\bprops[Bounds on Hessian sketch] 
\label{PropHessSketch}
For any convex set $\Constraint$ and any sketching matrix $\Sketch \in
\real^{\numproj \times \numobs}$, the Hessian sketch solution $\xhat$
satisfies the bound
\begin{align}
\label{EqnHessSketchBound}
\SEMI{\xhat - \xls} & \leq \frac{\ZSUP}{\ZINF} \; \SEMI{\xls}.
\end{align}
\eprops

For random sketching matrices, Proposition~\ref{PropHessSketch} can be
combined with probabilistic analysis to obtain high probability error
bounds. For a given tolerance parameter $\RHODEL \in (0,
\frac{1}{2}]$, consider the ``good event''
\begin{subequations}
\begin{align}
\label{EqnGoodEvent}
\Event(\RHODEL) & \defn \biggr \{ \ZINF \geq 1-\RHODEL, \mbox{ and }
\ZSUP \leq \frac{\RHODEL}{2} \biggr \}.
\end{align}
Conditioned on this event, Proposition~\ref{PropHessSketch} implies
that
\begin{align}
\label{EqnDeltaBound}
\SEMI{\xhat - \xls} & \leq \frac{\RHODEL}{2 \, (1 - \RHODEL)}
\SEMI{\xls} \; \leq \; \RHODEL \SEMI{\xls},
\end{align}
\end{subequations}
where the final inequality holds for all $\RHODEL \in (0,1/2]$.  

Thus, for a given family of random sketch matrices, we need to choose
the projection dimension $\numproj$ so as to ensure the event
$\Event{\RHODEL}$ holds for some $\RHODEL$.  For future reference, let
us state some known results for the cases of sub-Gaussian and ROS
sketching matrices.  We use $(c_0, c_1, c_2)$ to refer to numerical
constants, and we let $D = \dim(\Constraint)$ denote the dimension of
the space $\Constraint$.  In particular, we have $D = \usedim$ for
vector-valued estimation, and $D = \usedima \usedimb$ for matrix
problems.

Our bounds involve the ``size'' of the cone $\KCONELS$ previously
defined~\eqref{EqnTransformedCone}, as measured in terms of its
\emph{Gaussian width}
\begin{align}
\label{EqnPlainConeWidth}
\NonGaussComp(\KCONELS) & \defn \Exs_g \big[ \sup_{v \in \KCONELS \cap
    \Ball_2(1)} |\inprod{g}{v}| \big],
\end{align}
where $g \sim N(0, I_\numobs)$ is a standard Gaussian vector.  With
this notation, we have the following:
\blems[Sufficient conditions on sketch dimension~\cite{PilWai14a}]
\hfill
\label{LemSufficient}
\begin{enumerate}
\item[(a)] For sub-Gaussian sketch matrices, given a sketch size
  $\numproj > \frac{c_0}{\RHODEL^2} \NonGaussComp^2(\KCONELS)$, we have
\begin{subequations}
\begin{align}
\label{EqnGaussianSketch}
\mprob \big[ \Event(\RHODEL)] & \geq 1 - c_1 \CEXP{-c_2 \numproj
  \delta^2}.
\end{align}
\item[(b)] For randomized orthogonal system (ROS) sketches over the
  class of self-bounding cones, given a sketch size $\numproj >
  \frac{c_0 \, \log^4(D)}{\RHODEL^2} \NonGaussComp^2(\KCONELS)$, we
  have
\begin{align}
\mprob \big[ \Event(\RHODEL)] & \geq 1 - c_1 \CEXP{-c_2 \frac{\numproj
    \RHODEL^2}{\log^4(D)}}.
\end{align}
\end{subequations}
\end{enumerate}
\elems
The class of self-bounding cones is described more precisely in Lemma
8 of our earlier paper~\cite{PilWai14a}.  It includes among other
special cases the cones generated by unconstrained least-squares
(Example~\ref{ExaOrdinaryLeastSquares}), $\ell_1$-constrained least
squares (Example~\ref{ExaSparseLinear}), and least squares with
nuclear norm constraints (Example~\ref{ExaMultiTask}).  For these
cones, given a sketch size $\numproj > \frac{c_0 \,
  \log^4(D)}{\RHODEL^2} \NonGaussComp^2(\KCONELS)$, the Hessian sketch
applied with ROS matrices is guaranteed to return an estimate $\xhat$
such that $\SEMI{\xhat - \xls} \leq \RHODEL \SEMI{\xls}$ with high
probability.  This bound is an analogue of our earlier
bound~\eqref{EqnClassicalBound} for the classical sketch with
$\sqrt{f(\xls)}$ replaced by $\SEMI{\xls}$.  For this reason, we see
that the Hessian sketch alone suffers from the same deficiency as the
classical sketch: namely, it will require a sketch size $\numproj
\asymp \numobs$ in order to mimic the $\order(\numobs^{-1})$ accuracy
of the least-squares solution.


\subsection{Iterative Hessian sketch}

Despite the deficiency of the Hessian sketch itself, it serves as the
building block for an novel scheme---known as the iterative Hessian
sketch---that can be used to match the accuracy of the least-squares
solution using a reasonable sketch dimension.  Let begin by describing
the underlying intuition.  As summarized by the
bound~\eqref{EqnDeltaBound}, conditioned on the good event
$\Event(\RHODEL)$, the Hessian sketch returns an estimate with error
within a $\RHODEL$-factor of $\SEMI{\xls}$, where $\xls$ is the
solution to the original unsketched problem.  As show by
Lemma~\ref{LemSufficient}, as long as the projection dimension
$\numproj$ is sufficiently large, we can ensure that $\Event(\RHODEL)$
holds for some $\RHODEL \in (0, 1/2)$ with high probability.
Accordingly, given the current iterate $\xit{t}$, \emph{suppose that
  we can construct a new least-squares problem} for which the optimal
solution is $\xls - \xit{t}$.  Applying the Hessian sketch to this
problem will then produce a new iterate $\xit{t+1}$ whose distance to
$\xls$ has been reduced by a factor of $\RHODEL$.  Repeating this
procedure $\Tcrit$ times will reduce the initial approximation error
by a factor $\RHODEL^\Tcrit$.

With this intuition in place, we now turn a precise formulation of the
\emph{iterative Hessian sketch}.  Consider the optimization problem
\begin{align}
\label{EqnUproblem}
\uhat & = \arg \min_{u \in \Constraint - \xit{t}} \ENCMIN{ \frac{1}{2}
  \|\Amat u\|_2^2 - \inprod{\Amat^T (y - \Amat \xit{t})}{u}},
\end{align}
where $\xit{t}$ is the iterate at step $t$.  By construction, the
optimum to this problem is given by $\uhat = \xls - \xit{t}$.  We then
apply to Hessian sketch to this optimization
problem~\eqref{EqnUproblem} in order to obtain an approximation
$\xit{t+1} = \xit{t} + \uhat$ to the original least-squares solution
$\xls$ that is more accurate than $\xit{t}$ by a factor $\RHODEL \in
(0,1/2)$.  Recursing this procedure yields a sequence of iterates
whose error decays geometrically in $\RHODEL$. \\

\vspace*{.05in}

\noindent Formally, the iterative Hessian sketch algorithm takes the
following form:\\

\vspace*{.05in}

\begin{center}
\framebox[.99\textwidth]{
\parbox{.95\textwidth}{ 
\noindent {\bf{Iterative Hessian sketch (IHS):}}  Given an iteration
number $\Tcrit \geq 1$:
\begin{enumerate}
\item[(1)] Initialize at $\xit{0} = 0$.
\item[(2)] For iterations $t = 0, 1, 2, \ldots, \Tcrit-1$, generate an
  independent sketch matrix \mbox{$\SketchIt{t+1} \in \real^{\numproj
      \times \numobs}$,} and perform the update
\begin{align}
\label{EqnXitUpdate}
\xit{t+1} & = \arg \min_{x \in \Constraint} \ENCMIN{ 
\frac{1}{2
    \numproj} \|\SketchIt{t+1} \Amat(x - \xit{t})\|_2^2 -
  \inprod{\Amat^T (y - \Amat \xit{t})}{x}}.
\end{align}
\item[(3)]  Return the estimate $\xhat = \xit{\Tcrit}$.
\end{enumerate}
}
}
\end{center}
\vspace*{.1in}

\noindent The following theorem summarizes the key properties of this
algorithm.  It involves the sequence $\{\ZINFIT{t}, \ZSUPIT{t}
\}_{t=1}^\Tcrit$, where the quantities $\ZINF$ and $\ZSUP$ were
previously defined in equations~\eqref{EqnDefnZinf}
and~\eqref{EqnDefnZsup}.  In addition, as a generalization of the
event~\eqref{EqnGoodEvent}, we define the sequence of ``good'' events
\begin{align}
\EVENTIT{t}(\RHODEL) & \defn \biggr \{ \ZINFIT{t} \geq 1-\RHODEL,
\mbox{ and } \ZSUPIT{t} \leq \frac{\RHODEL}{2} \biggr \} \qquad
\mbox{for $t = 1, \ldots, \Tcrit$.}
\end{align}
With this notation, we have the following guarantee:
\btheos[Guarantees for iterative Hessian sketch]
\label{ThmOptIterative}
The final solution $\xhat = \xit{\Tcrit}$ satisfies the bound
\begin{subequations}
\begin{align}
\label{EqnProdBound}
\SEMI{\xhat - \xls} & \leq \Big \{ \prod_{t=1}^{\Tcrit}
\frac{\ZSUPIT{t}}{\ZINFIT{t}} \Big \} \; \SEMI{\xls}.
\end{align}
Consequently, conditioned on the event $\cap_{t=1}^\Tcrit
\EVENTIT{t}(\RHODEL)$ for some $\RHODEL \in (0, 1/2)$, we have
\begin{align}
\label{EqnFinalBound}
\SEMI{\xhat - \xls} & \leq \RHODEL^\Tcrit \; \SEMI{\xls}.
\end{align}
\end{subequations}
\etheos
\noindent Note that for any $\RHODEL \in (0, 1/2)$, then event
$\EVENTIT{t}(\RHODEL)$ implies that $\frac{\ZSUPIT{t}}{\ZINFIT{t}}
\leq \RHODEL$, so that the bound~\eqref{EqnFinalBound} is an immediate
consequence of the product bound~\eqref{EqnProdBound}.\\


Lemma~\ref{LemSufficient} can be combined with the union bound in
order to ensure that the compound event $\cap_{t=1}^\Tcrit
\EVENTIT{t}(\RHODEL)$ holds with high probability over a sequence of
$N$ iterates, as long as the sketch size is lower bounded as $\numproj
\geq \frac{c_0}{\RHODEL^2} \Width^2(\KCONELS) \log^4(D) + \log
\Tcrit$.  Based on the bound~\eqref{EqnFinalBound}, we then expect to
observe geometric convergence of the iterates. 

\begin{figure}[h]
\begin{center}
\begin{tabular}{cc}
\widgraph{.5\textwidth}{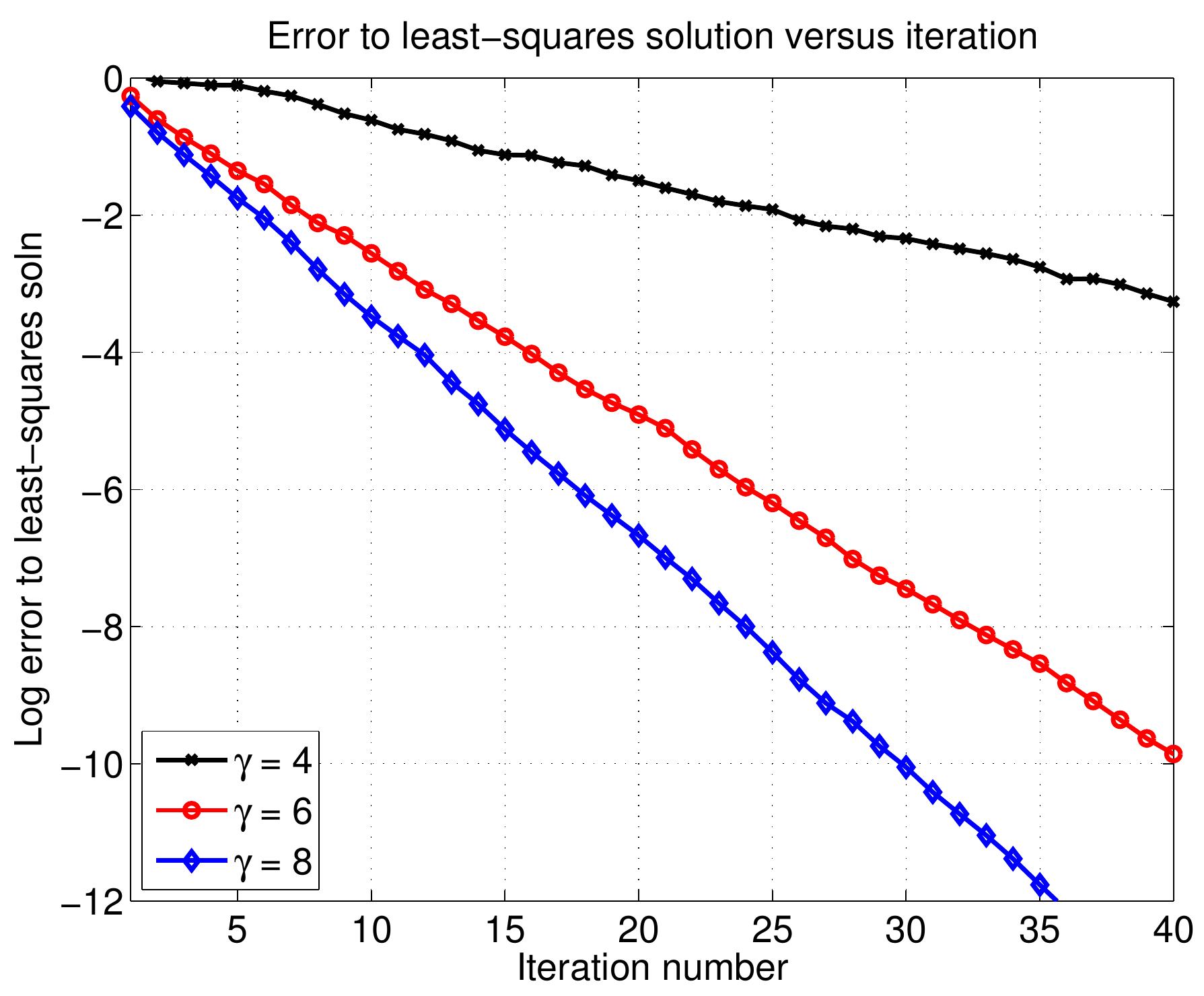} &
\widgraph{.5\textwidth}{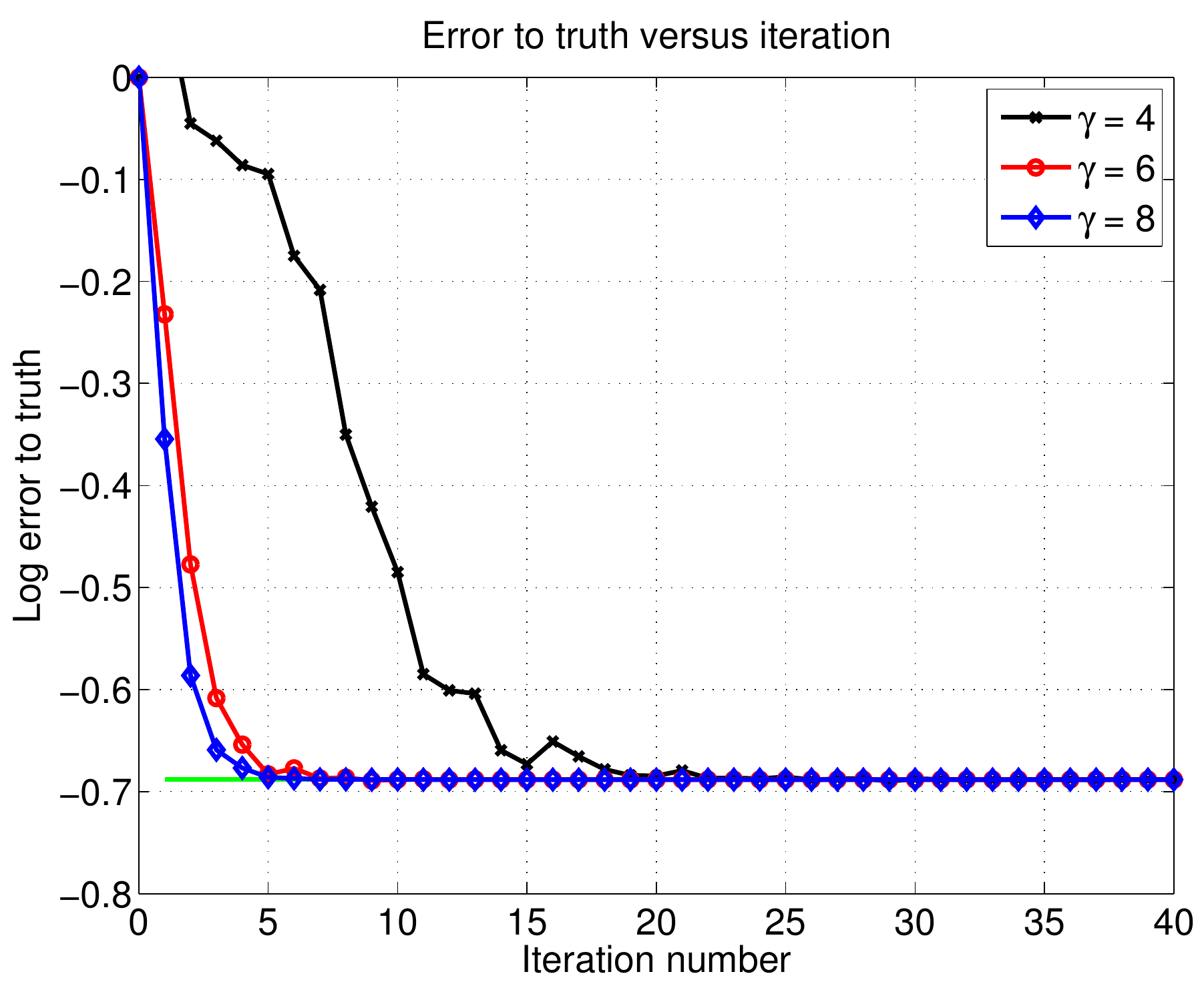} \\
(a) & (b)
\end{tabular}
\caption{Simulations of the IHS algorithm for an unconstrained
  least-squares problem with noise variance $\sigma^2 = 1$, and of
  dimensions $(\usedim, \numobs) = (200, 6000)$.  Simulations based on
  sketch sizes $\numproj = \gamma \usedim$, for a parameter $\gamma >
  0$ to be set.  (a) Plots of the log error $\SEMI{\xit{t} - \xls}$
  versus the iteration number $t$.  Three different curves for $\gamma
  \in \{4, 6, 8 \}$.  Consistent with the theory, the convergence is
  geometric, with the rate increasing as the sampling factor $\gamma$
  is increased. (b) Plots of the log error $\SEMI{\xit{t} - \xstar}$
  versus the iteration number $t$.  Three different curves for $\gamma
  \in \{4, 6, 8 \}$.  As expected, all three curves flatten out at the
  level of the least-squares error $\SEMI{\xls - \xstar} = 0.20
  \approx \sqrt{\sigma^2 \usedim/\numobs}$.}
\label{FigOptLS}
\end{center}
\end{figure}

In order to test this prediction, we implemented the IHS algorithm
using Gaussian sketch matrices, and applied it to an unconstrained
least-squares problem based on a data matrix with dimensions
\mbox{$(\usedim, \numobs) = (200, 6000)$} and noise variance $\sigma^2
= 1$. As shown in Appendix~\ref{AppExamples}, the Gaussian width of
$\KCONELS$ is proportional to $\usedim$, so that
Lemma~\ref{LemSufficient} shows that it suffices to choose a
projection dimension $\numproj \succsim \gamma \usedim$ for a
sufficiently large constant $\gamma$.  Panel (a) of
Figure~\ref{FigOptLS} illustrates the resulting convergence rate of
the IHS algorithm, measured in terms of the error $\SEMI{\xit{t} -
  \xls}$, for different values $\gamma \in \{4, 6, 8\}$. As predicted
by Theorem~\ref{ThmOptIterative}, the convergence rate is geometric
(linear on the log scale shown), with the rate increasing as the
parameter $\gamma$ is increased.

\vspace*{.1in}

Assuming that the sketch dimension has been chosen to ensure geometric
convergence, Theorem~\ref{ThmOptIterative} allows us to specify, for a
given target accuracy $\varepsilon \in (0,1)$, the number of
iterations required.
\bcors
\label{CorGeneral}
Fix some $\RHODEL \in (0, 1/2)$, and choose a sketch dimension
$\numproj > \frac{c_0 \log^4(D)}{\RHODEL^2}
\NonGaussComp^2(\KCONELS)$.  If we apply the IHS algorithm for
$\Tcrit(\RHODEL, \varepsilon) \defn 1 +
\frac{\log(1/\varepsilon)}{\log(1/\RHODEL)}$ steps, then the output
$\xhat = \xit{\Tcrit}$ satisfies the bound
\begin{align}
\frac{\SEMI{\xhat - \xls}}{\SEMI{\xls}} & \leq \varepsilon
\end{align}
with probability at least $1 - \plaincon_1 \Tcrit(\RHODEL,
\varepsilon) \CEXP{-\plaincon_2 \frac{\numproj
    \RHODEL^2}{\log^4(D)}}$.
\ecors

\noindent This corollary is an immediate consequence of
Theorem~\ref{ThmOptIterative} combined with Lemma~\ref{LemSufficient},
and it holds for both ROS and sub-Gaussian sketches.  (In the latter
case, the additional $\log(D)$ terms may be omitted.)  Combined with
bounds on the width function $\NonGaussComp(\KCONELS)$, it leads to a
number of concrete consequences for different statistical models, as
we illustrate in the following section.

One way to understand the improvement of the IHS algorithm over the
classical sketch is as follows.  Fix some error tolerance $\varepsilon
\in (0,1)$.  Disregarding logarithmic factors, our previous
results~\cite{PilWai14a} on the classical sketch then imply that a
sketch size $\numproj \succsim \varepsilon^{-2} \:
\NonGaussComp^2(\KCONELS)$ is sufficient to produce a
$\varepsilon$-accurate solution approximation.  In contrast,
Corollary~\ref{CorGeneral} guarantees that a sketch size $\numproj
\succsim \log(1/\varepsilon) \; \NonGaussComp^2(\KCONELS)$ is
sufficient.  Thus, the benefit is the reduction from
$\varepsilon^{-2}$ to $\log(1/\varepsilon)$ scaling of the required
sketch size.

It is worth noting that in the absence of constraints, the
least-squares problem reduces to solving a linear system, so that
alternative approaches are available.  For instance, one can use a
randomized sketch to obtain a preconditioner, which can then be used
within the conjugate gradient method.  As shown in past
work~\cite{Rokhlin2008,Avron2010}, two-step methods of this type can
lead to same reduction of $\varepsilon^{-2}$ dependence to
$\log(1/\varepsilon)$.  However, a method of this type is very
specific to unconstrained least-squares, whereas the procedure
described in this paper is generally applicable to least-squares over
any compact, convex constraint set.


\subsection{Computational and space complexity}

Let us now make a few comments about the computational and space
complexity of implementing the IHS algorithm using ROS sketches (e.g.,
such as those based on the fast Hadamard transform).  For a given
sketch size $\numproj$, at iteration $t$, the IHS algorithm requires
$\order(\numobs \usedim \log(\numproj))$ basic operation for computing
the data sketch $\Sketch ^{t+1} \Amat$ and also $\order(\numobs
\usedim)$ operations to compute $\Amat^T (y - \Amat x^t)$.
Consequently, if we run the algorithm for $\Tcrit$ iterations, then
the overall complexity is $\order \big( (\numobs \usedim
\log(\numproj) + \COMP(\numproj, \usedim) ) \, \Tcrit\big)$, where
$\COMP(\numproj, \usedim)$ is the complexity of solving the $m\times
d$ dimensional problem in the update~\eqref{EqnXitUpdate}. The total
space used scales as $\order(\numproj \usedim)$.

If we want to obtain estimates with accuracy $\varepsilon$, then we
need to perform $N \asymp \log (1/\varepsilon)$ iterations in total.
Moreover, for ROS sketches, we need to choose $\numproj \succsim
\Width^2(\KCONELS) \log^4(d)$.  Consequently, it only remains to bound
the Gaussian width $\Width$ in order to specify complexities that
depend only on the pair $(\numobs, \usedim)$.

\paragraph{Unconstrained least-squares:}
For an unconstrained problem with $\numobs > \usedim$, the Gaussian
width can be bounded as $\Width^2(\KCONELS) \precsim \usedim$, and the
complexity of the solving the sub-problem~\eqref{EqnXitUpdate} can be
bounded as $\usedim^3$. Thus, the overall complexity of computing an
$\varepsilon$-accurate solution scales as $\order(\numobs \usedim
\log(\usedim) + \usedim^3) \log(1/\varepsilon)$, and the space
required is $\order(\usedim^2)$.

\paragraph{Sparse least-squares:}   As will be shown in 
Section~\ref{SecConsequencesSparse}, in certain cases, the cone
$\KCONELS$ can have substantially lower complexity than the
unconstrained case.  For instance, if the solution is sparse, say with
$\spindex$ non-zero entries and the least-squares program involves an
$\ell_1$-constraint, then we have $\Width^2(\KCONELS) \precsim
\spindex \log \usedim$.  Using a standard interior point method to
solve the sketched problem, the total complexity for obtaining an
$\varepsilon$-accurate solution is upper bounded by $\order((\numobs \usedim \log(\spindex ) + \spindex^2
\usedim \log^2(\usedim) )\log(1/\varepsilon))$, along with an
$\order(\spindex \usedim \log(\usedim))$ space complexity. The sparsity $\spindex$ is not known a priori, however we note that there exists efficient bounds on $\spindex$ which can be computed in $\order(nd)$ time (see e.g. \cite{el2011safe}). Another approach is to impose some conditions on the design matrix $\Amat$ which will guarantee support recovery, i.e.,  number of nonzero entries of $x^*$ equals $\spindex$, (e.g., see  \cite{Wainwright06}).


\section{Consequences for concrete models}
\label{SecConsequences}

In this section, we derive some consequences of
Corollary~\ref{CorGeneral} for particular classes of least-squares
problems.  Our goal is to provide empirical confirmation of the
sharpness of our theoretical predictions, namely the minimal sketch
dimension required in order to match the accuracy of the original
least-squares solution.

\subsection{Unconstrained least squares}
\label{SecConsequencesUnc}

We begin with the simplest case, namely the unconstrained
least-squares problem ($\Constraint = \real^\usedim$).  For a given
pair $(\numobs, \usedim)$ with $\numobs > \usedim$, we generated a
random ensemble of least-square problems according to the following
procedure:
\begin{itemize}
\item first, generate a random data matrix $\Amat \in \real^{\numobs
  \times \usedim}$ with i.i.d. $N(0,1)$ entries
\item second, choose a regression vector $\xstar$ uniformly at random
  from the sphere $\Sphere{\usedim}$
\item third, form the response vector $y = \Amat \xstar + w$, where $w
  \sim N(0, \sigma^2 I_\numobs)$ is observation noise with $\sigma =
  1$.  
\end{itemize}
As discussed following Lemma~\ref{LemSufficient}, for this class of
problems, taking a sketch dimension \mbox{$\numproj \succsim
  \frac{\usedim}{\RHODEL^2}$} guarantees $\RHODEL$-contractivity of
the IHS iterates with high probability.  Consequently, we can obtain a
$\varepsilon$-accurate approximation to the original least-squares
solution by running roughly $\log(1/\varepsilon)/\log(1/\RHODEL)$
iterations.

Now how should the tolerance $\varepsilon$ be chosen?  Recall that the
underlying reason for solving the least-squares problem is to
approximate $\xstar$.  Given this goal, it is natural to measure the
approximation quality in terms of $\SEMI{\xit{t} - \xstar}$.  Panel
(b) of Figure~\ref{FigOptLS} shows the convergence of the iterates to
$\xstar$.  As would be expected, this measure of error levels off at
the ordinary least-squares error
\begin{align*}
\SEMI{\xls - \xstar}^2 \asymp \frac{\sigma^2 \usedim}{\numobs} \approx
0.10.
\end{align*}
Consequently, it is reasonable to set the tolerance parameter
proportional to $\sigma^2 \frac{\usedim}{\numobs}$, and then perform
roughly $1 + \frac{\log(1/\varepsilon)}{\log(1/\RHODEL)}$ steps.  The
following corollary summarizes the properties of the resulting
procedure:
\bcors
\label{CorOptLS}
For some given $\RHODEL \in (0, 1/2)$, suppose that we run the IHS
algorithm for \mbox{$\Tcrit = 1 + \lceil \frac{\log \sqrt{\numobs} \;
    \frac{\SEMI{\xls}}{\sigma}} {\log(1/\RHODEL)} \rceil$} iterations
using $\numproj = \frac{c_0}{\RHODEL^2} \usedim$ projections per
round. Then the output $\xhat$ satisfies the bounds
\begin{align}
\label{EqnXlserrLS}
\SEMI{\xhat - \xls} \leq \sqrt{\frac{\sigma^2 \usedim}{\numobs}},
\qquad \mbox{and} \qquad \SEMI{\xit{\Tcrit} - \xstar} \leq
\sqrt{\frac{\sigma^2 \usedim}{\numobs}} + \SEMI{\xls - \xstar}
\end{align}
with probability greater than $1 - \plaincon_1 \, \Tcrit \,
\CEXP{-\plaincon_2 \frac{\numproj \RHODEL^2}{\log^4(\usedim)}}$.
\ecors

\begin{figure}[h]
\begin{center}
\widgraph{.6\textwidth}{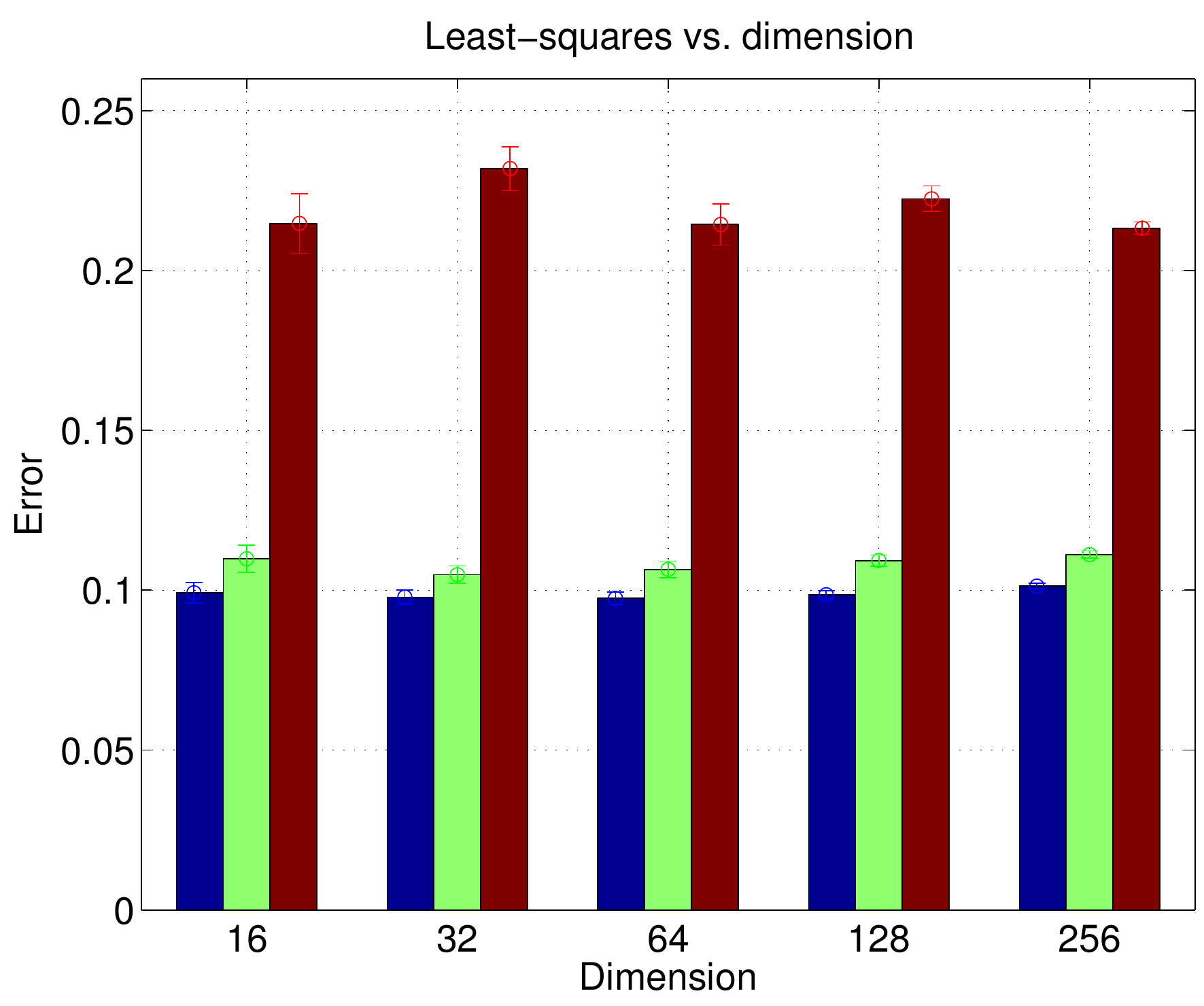}
\caption{Simulations of the IHS algorithm for unconstrained
  least-squares.  In these experiments, we generated random
  least-squares problem of dimensions $\usedim \in \{16, 32, 64, 128,
  256\}$, on all occasions with a fixed sample size $\numobs = 100
  \usedim$. The initial least-squares solution has error $\SEMI{\xls -
    \xstar} \approx 0.10$, as shown by the blue bars.  We then ran the
  IHS algorithm for $\Tcrit = 4$ iterations with a sketch size
  $\numproj = 6 \usedim$. As shown by the green bars, these sketched
  solutions show an error $\SEMI{\xhat - \xstar}\approx 0.11$
  independently of dimension, consistent with the predictions of
  Corollary~\ref{CorOptLS}.  Finally, the red bars show the error in
  the classical sketch, based on a sketch size $M = \Tcrit \numproj =
  24 \usedim$, corresponding to the total number of projections used
  in the iterative algorithm.  This error is roughly twice as large.}
\label{FigScaleLS}
\end{center}
\end{figure}

In order to confirm the predicted bound~\eqref{EqnXlserrLS} on the
error $\SEMI{\xhat - \xls}$, we performed a second experiment.  Fixing
$\numobs = 100 \usedim$, we generated $T = 20$ random least squares
problems from the ensemble described above with dimension $\usedim$
ranging over $\{32, 64, 128, 256, 512 \}$.  By our previous choices,
the least-squares estimate should have error $\|\xls - \xstar\|_2
\approx \sqrt{\frac{\sigma^2 \usedim}{\numobs}} = 0.1$ with high
probability, independently of the dimension $\usedim$.  This predicted
behavior is confirmed by the blue bars in Figure~\ref{FigScaleLS}; the
bar height corresponds to the average over $T = 20$ trials, with the
standard errors also marked.  On these same problem instances, we also
ran the IHS algorithm using $\numproj = 6 \usedim$ samples per
iteration, and for a total of
\begin{align*}
\Tcrit & = 1 + \lceil \frac{\log \big( \sqrt{\frac{\numobs}{\usedim}}
  \big)}{\log 2} \rceil \; =  4 \qquad \mbox{iterations.}
\end{align*}
Since $\SEMI{\xls - \xstar} \asymp \sqrt{\frac{\sigma^2
    \usedim}{\numobs}} \approx 0.10$, Corollary~\ref{CorOptLS} implies
that with high probability, the sketched solution $\xhat =
\xit{\Tcrit}$ satisfies the error bound
\begin{align*}
\|\xhat - \xstar\|_2 & \leq \plaincon'_0 \; \sqrt{\frac{\sigma^2
    \usedim}{\numobs}}
\end{align*}
for some constant $\plaincon'_0 > 0$.  This prediction is confirmed by
the green bars in Figure~\ref{FigScaleLS}, showing that $\SEMI{\xhat -
  \xstar} \approx 0.11$ across all dimensions.  Finally, the red bars
show the results of running the classical sketch with a sketch
dimension of $(6 \times 4) \usedim = 24 \usedim$ sketches,
corresponding to the total number of sketches used by the IHS
algorithm.  Note that the error is roughly twice as large.

\subsection{Sparse least-squares}
\label{SecConsequencesSparse}
We now turn to a study of an $\ell_1$-constrained form of
least-squares, referred to as the Lasso or relaxed basis
pursuit~\cite{Chen98,Tibshirani96}.  In particular, consider the
convex program
\begin{align}
\xls & = \arg \min_{\|x\|_1 \leq R} \big \{ \frac{1}{2} \|\yvec -
\Amat x\|_2^2 \big \},
\end{align}
where $R > 0$ is a user-defined radius.  This estimator is well-suited
to the problem of sparse linear regression, based on the observation
model $y = \Amat \xstar + w$, where $\xstar$ has at most $\kdim$
non-zero entries, and $\Amat \in \real^{\numobs \times \usedim}$ has
i.i.d. $N(0,1)$ entries.  For the purposes of this illustration, we
assume\footnote{In practice, this unrealistic assumption of exactly
  knowing $\|\xstar\|_1$ is avoided by instead considering the
  $\ell_1$-penalized form of least-squares, but we focus on the
  constrained case to keep this illustration as simple as possible.}
that the radius is chosen such that $R = \|\xstar\|_1$.

\begin{figure}[t]
\begin{center}
\begin{tabular}{cc}
\widgraph{.5\textwidth}{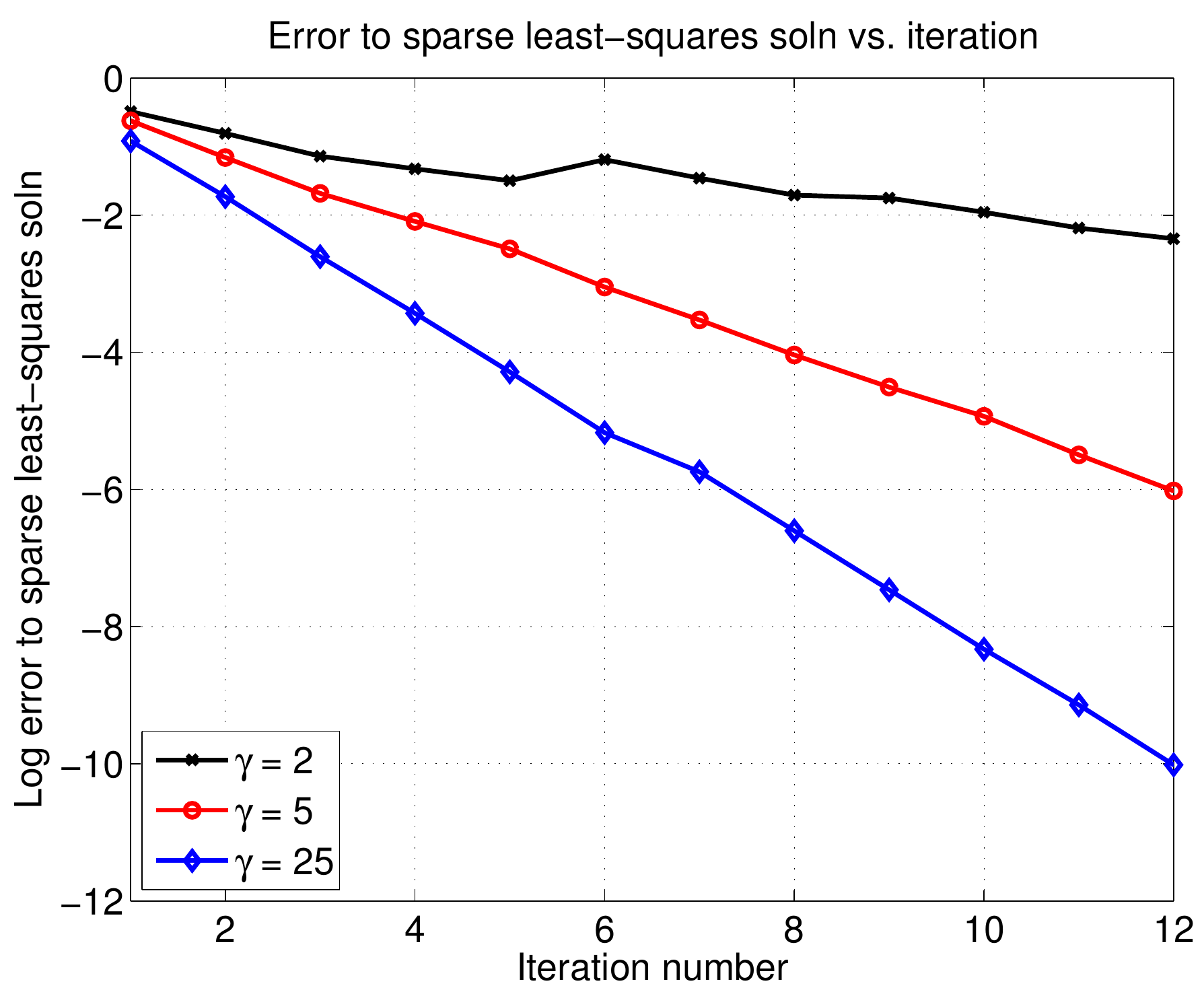} &
\widgraph{.5\textwidth}{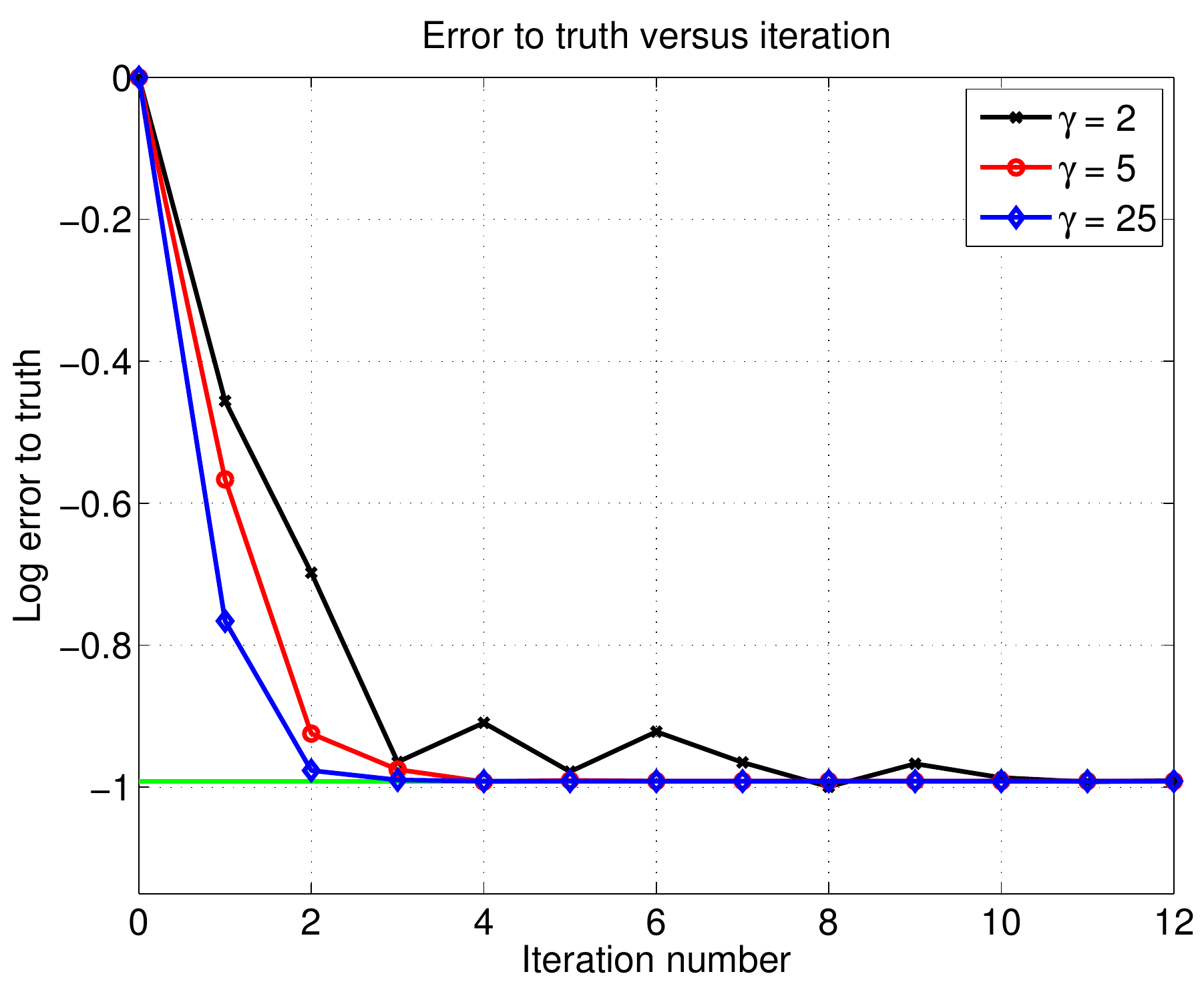} \\
(a) & (b)
\end{tabular}
\caption{Simulations of the IHS algorithm for a sparse least-squares
  problem with noise variance $\sigma^2 = 1$, and of dimensions
  $(\usedim, \numobs, \kdim) = (256, 8872, 32)$.  Simulations based on
  sketch sizes $\numproj = \gamma \kdim \log \usedim$, for a parameter
  $\gamma > 0$ to be set.  (a) Plots of the log error $\|\xit{t} -
  \xls\|_2$ versus the iteration number $t$.  Three different curves
  for $\gamma \in \{2, 5, 25 \}$.  Consistent with the theory, the
  convergence is geometric, with the rate increasing as the sampling
  factor $\gamma$ is increased. (b) Plots of the log error $\|\xit{t}
  - \xstar\|_2$ versus the iteration number $t$.  Three different
  curves for $\gamma \in \{2, 5, 25 \}$.  As expected, all three
  curves flatten out at the level of the least-squares error $\|\xls -
  \xstar\|_2 = 0.10 \approx \sqrt{\frac{\kdim \log (e
      \usedim/\kdim)}{\numobs}}$.}
\label{FigOptSparse}
\end{center}
\end{figure}

Under these conditions, the proof of Corollary~\ref{CorOptSparse}
shows that a sketch size $\numproj \geq \gamma \, \kdim \log
\big(\frac{e \usedim}{\kdim} \big)$ suffices to guarantee geometric
convergence of the IHS updates.  Panel (a) of
Figure~\ref{FigOptSparse} illustrates the accuracy of this prediction,
showing the resulting convergence rate of the the IHS algorithm,
measured in terms of the error $\SEMI{\xit{t} - \xls}$, for different
values $\gamma \in \{2, 5, 25\}$. As predicted by
Theorem~\ref{ThmOptIterative}, the convergence rate is geometric
(linear on the log scale shown), with the rate increasing as the
parameter $\gamma$ is increased.

As long as $\numobs \succsim \kdim \log \big(\frac{e \usedim}{\kdim}
\big)$, it also follows as a corollary of Proposition~\ref{PropMLE}
that
\begin{align}
\SEMI{\xls - \xstar}^2 \precsim \frac{\sigma^2 \kdim \log \big(\frac{e
    \usedim}{\kdim} \big)}{\numobs}.
\end{align}
with high probability.  This bound suggests an appropriate choice for
the tolerance parameter $\varepsilon$ in
Theorem~\ref{ThmOptIterative}, and leads us to the following
guarantee.

\bcors
\label{CorOptSparse}
For the stated random ensemble of sparse linear regression problems,
suppose that we run the IHS algorithm for $\Tcrit = 1 + \lceil
\frac{\log \sqrt{\numobs} \; \frac{\SEMI{\xls}}{\sigma}}
     {\log(1/\RHODEL)} \rceil $ iterations using $\numproj =
     \frac{c_0}{\RHODEL^2} \kdim \log \big(\frac{e \usedim}{\kdim}
     \big)$ projections per round. Then with probability greater than
     $1 - \plaincon_1 \, \Tcrit \, \CEXP{-\plaincon_2 \frac{\numproj
         \RHODEL^2}{\log^4(\usedim)}}$, the output $\xhat$ satisfies
     the bounds
\begin{align}
\label{EqnOptSparseBound}
\SEMI{\xhat - \xls} \leq \sqrt{\frac{\sigma^2 \kdim \log \big(\frac{e
      \usedim}{\kdim} \big)}{\numobs}} \qquad \mbox{and} \qquad
\SEMI{\xit{\Tcrit} - \xstar} \leq \sqrt{\frac{\sigma^2 \kdim \log
    \big(\frac{e \usedim}{\kdim} \big)}{\numobs}} + \SEMI{\xls -
  \xstar}.
\end{align}
\ecors

\begin{figure}[h]
\begin{center}
\widgraph{.6\textwidth}{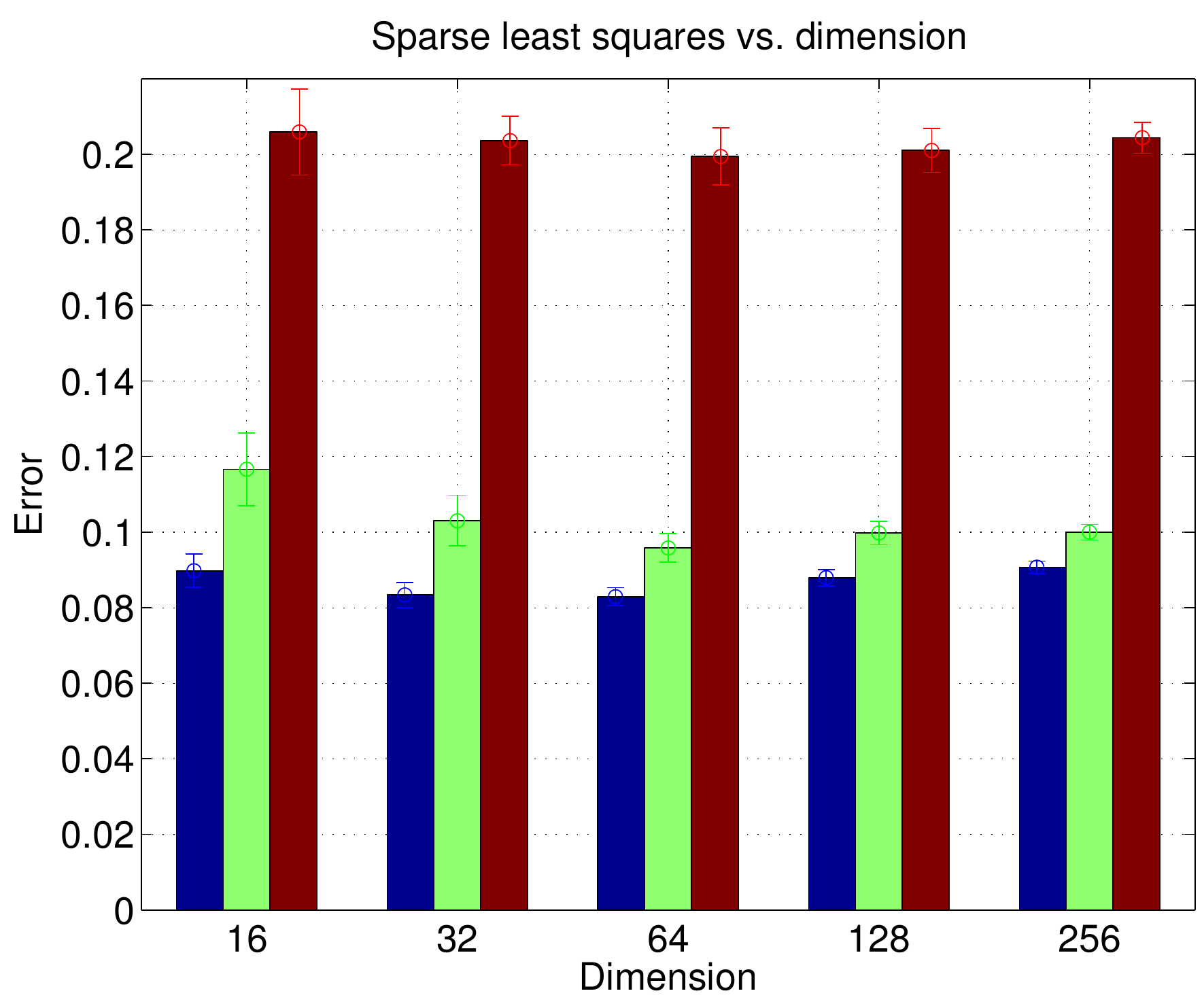}
\caption{Simulations of the IHS algorithm for $\ell_1$-constrained
  least-squares.  In these experiments, we generated random sparse
  least-squares problem of dimensions \mbox{$\usedim \in \{16, 32, 64,
    128, 256\}$} and sparsity $\kdim = \lceil 2 \sqrt{\usedim}
  \rceil$, on all occasions with a fixed sample size $\numobs = 100
  \kdim \log \big( \frac{e \usedim}{\kdim} \big)$. The initial Lasso
  solution has error $\|\xls - \xstar\|_2 \approx 0.10$, as shown by
  the blue bars.  We then ran the IHS algorithm for $\Tcrit = 4$
  iterations with a sketch size $\numproj = 4 \kdim \log \big(\frac{e
    \usedim}{\kdim} \big)$. These sketched solutions show an error
  $\SEMI{\xhat - \xstar} \approx 0.11$ independently of dimension,
  consistent with the predictions of Corollary~\ref{CorOptSparse}.
  Red bars show the error in the naive sketch estimate, using a sketch
  of size $M = \Tcrit \numproj = 16 \kdim \log \big ( \frac{ e
    \usedim}{\kdim} \big)$, equal to the total number of random
  projections used by the IHS algorithm.  The resulting error is
  roughly twice as large. }
\label{FigScaleSparse}
\end{center}
\end{figure}

In order to verify the predicted bound~\eqref{EqnOptSparseBound} on
the error $\SEMI{\xhat - \xls}$, we performed a second experiment.
Fixing $\numobs = 100 \kdim \log \big( \frac{e \usedim}{\kdim} \big)$.
we generated $T = 20$ random least squares problems (as described
above) with the regression dimension ranging as $\usedim \in \{32, 64,
128, 256 \}$, and sparsity $\kdim = \lceil 2 \sqrt{\usedim} \rceil$.
Based on these choices, the least-squares estimate should have error
$\SEMI{\xls - \xstar} \approx \sqrt{\frac{\sigma^2 \kdim \log
    \big(\frac{e \usedim}{\kdim} \big)}{\numobs}} = 0.1$ with high
probability, independently of the pair $(\kdim, \usedim)$.  This
predicted behavior is confirmed by the blue bars in
Figure~\ref{FigScaleSparse}; the bar height corresponds to the average
over $T = 20$ trials, with the standard errors also marked.

On these same problem instances, we also ran the IHS algorithm using
$\Tcrit = 4$ iterations with a sketch size $\numproj = 4 \kdim \log
\big(\frac{e \usedim}{\kdim} \big)$.  Together with our earlier
calculation of $\SEMI{\xls - \xstar}$, Corollary~\ref{CorOptLS}
implies that with high probability, the sketched solution $\xhat =
\xit{\Tcrit}$ satisfies the error bound
\begin{align*}
\SEMI{\xhat - \xstar} & \leq \plaincon_0 \; \sqrt{\frac{\sigma^2 \kdim
    \log \big( \frac{e \usedim}{\kdim} \big)}{\numobs}}
\end{align*}
for some constant $\plaincon_0 \in (1,2]$.  This prediction is
  confirmed by the green bars in Figure~\ref{FigScaleSparse}, showing
  that $\SEMI{\xhat - \xstar} \succsim 0.11$ across all dimensions.
  Finally, the green bars in Figure~\ref{FigScaleSparse} show the
  error based on using the naive sketch estimate with a total of $M =
  \Tcrit \numproj$ random projections in total; as with the case of
  ordinary least-squares, the resulting error is roughly twice as
  large. We also note that a similar bound also applies to problems where a parameter constrained to unit simplex is estimated, e.g., in portfolio analysis and density estimation \cite{Markowitz59,pilanci2012recovery}.


\subsection{Matrix estimation with nuclear norm constraints}

We now turn to the study of nuclear-norm constrained form of
least-squares matrix regression.  This class of problems has proven
useful in many different application areas, among them matrix
completion, collaborative filtering, multi-task learning and control
theory
(e.g.,~\cite{Fa02,YuaEkiLuMon07,Bach08b,RecFazPar10,NegWai10b}).  In
particular, let us consider the convex program
\begin{align}
\label{EqnNucNormLS}
\Xls & = \arg \min_{X \in \real^{\usedima \times \usedimb}}
\ENCMIN{\frac{1}{2} \fronorm{Y - \Amat X}^2} \qquad \mbox{such that
  $\nucnorm{X} \leq R$,}
\end{align}
where $R > 0$ is a user-defined radius as a regularization parameter.
 

\subsubsection{Simulated data}

Recall the linear observation model previously introduced in
Example~\ref{ExaMultiTask}: we observe the pair $(Y, \Amat)$ linked
according to the linear $Y = \Amat \Xstar + W$, where the unknown
matrix $\Xstar\in \real^{\usedima \times \usedimb}$ is an unknown
matrix of rank $\rdim$.  The matrix $W$ is observation noise, formed
with i.i.d. $N(0, \sigma^2)$ entries.  This model is a special case of
the more general class of matrix regression problems~\cite{NegWai10b}.
As shown in Appendix~\ref{AppExamples}, if we solve the nuclear-norm
constrained problem with $R = \nucnorm{X^*}$, then it produces a
solution such that $\Exs \big[ \fronorm{\Xls - \Xstar}^2] \precsim
\sigma^2 \frac{\rdim \, (\usedima + \usedimb)}{\numobs}$.  The
following corollary characterizes the sketch dimension and iteration
number required for the IHS algorithm to match this scaling up to a
constant factor.
\bcors[IHS for nuclear-norm constrained least squares]
\label{CorLowRank}
Suppose that we run the IHS algorithm for $\Tcrit = 1 + \lceil
\frac{\log \sqrt{\numobs} \; \frac{\SEMI{\Xls}}{\sigma}}
     {\log(1/\RHODEL)} \rceil $ iterations using $\numproj =
     {c_0}{\RHODEL^2} \rdim \big(\usedima + \usedimb \big)$
     projections per round. Then with probability greater than $1 -
     \plaincon_1 \, \Tcrit \, \CEXP{-\plaincon_2 \frac{\numproj
         \RHODEL^2}{\log^4(\usedima \usedimb)}}$, the output
     $X^\Tcrit$ satisfies the bound
\begin{align}
\label{EqnOptLowRankBound}
\SEMI{X^{\Tcrit} - \Xstar} \leq \sqrt{\frac{\sigma^2 \rdim
    \big(\usedima + \usedimb \big)}{\numobs}} + \SEMI{\Xls - \Xstar}.
\end{align}
\ecors
\noindent

We have also performed simulations for low-rank matrix estimation, and
observed that the IHS algorithm exhibits convergence behavior
qualitatively similar to that shown in Figures~\ref{FigScaleLS}
and~\ref{FigScaleSparse}.  Similarly, panel (a) of
Figure~\ref{FigLowRank} compares the performance of the IHS and
classical methods for sketching the optimal solution over a range of
row sizes $\numobs$.  As with the unconstrained least-squares results
from Figure~\ref{FigLeastSquaresFixedM}, the classical sketch is very
poor compared to the original solution whereas the IHS algorithm
exhibits near optimal performance.


\subsubsection{Application to multi-task learning}

To conclude, let us illustrate the use of the IHS algorithm in
speeding up the training of a classifier for facial expressions.  In
particular, suppose that our goal is to separate a collection of
facial images into different groups, corresponding either to distinct
individuals or to different facial expressions.  One approach would be
to learn a different linear classifier ($a \mapsto \inprod{a}{x}$) for
each separate task, but since the classification problems are so
closely related, the optimal classifiers are likely to share
structure.  One way of capturing this shared structure is by
concatenating all the different linear classifiers into a matrix, and
then estimating this matrix in conjunction with a nuclear norm
penalty~\cite{Amit07,Argyriou08}.

\begin{figure}[h]
\begin{center}
\begin{tabular}{cc}
\includegraphics[trim = 0mm 0mm 0mm 1mm, clip, width=0.7\textwidth]{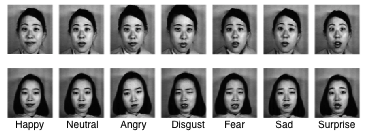}
\end{tabular}
\caption{Japanese Female Facial Expression (JAFFE) Database: The JAFFE
  database consists of 213 images of 7 different emotional facial
  expressions (6 basic facial expressions + 1 neutral) posed by 10
  Japanese female models.}
\label{FigJAFFE}
\end{center}
\end{figure}

In more detail, we performed a simulation study using the The Japanese
Female Facial Expression (JAFFE) database~\cite{Lyons98}. It consists
of $N = 213$ images of $7$ facial expressions ($6$ basic facial
expressions + $1$ neutral) posed by $10$ different Japanese female
models; see Figure~\ref{FigJAFFE} for a few example images.  We
performed an approximately $80:20$ split of the data set into $\ntrain
= 170$ training and $\ntest = 43$ test images respectively. Then we
consider classifying each facial expression and each female model as a
separate task which gives a total of $\dtask = 17$ tasks.  For each
task $j = 1, \ldots, \dtask$, we construct a linear classifier of the
form $a \mapsto \sign(\inprod{a}{x_j})$, where $a \in \real^\usedim$
denotes the vectorized image features given by Local Phase
Quantization~\cite{Ojansivu08}.  In our implementation, we fixed the
number of features $\usedim = 32$.  Given this set-up, we train the
classifiers in a joint manner, by optimizing simultaneously over the
matrix $X \in \real^{\usedim \times \dtask}$ with the classifier
vector $x_j \in \real^{\usedim}$ as its $j^{th}$ column.  The image
data is loaded into the matrix $A \in \real^{\ntrain \times \usedim}$,
with image feature vector $a_i \in \real^\usedim$ in column $i$ for $i
= 1, \ldots, \ntrain$.  Finally, the matrix $Y \in \{-1, +1\}^{\ntrain
  \times \dtask}$ encodes class labels for the different
  classification problems.  These instantiations of the pair $(Y, X)$
  give us an optimization problem of the form~\eqref{EqnNucNormLS},
  and we solve it over a range of regularization radii $R$.

\begin{figure}[h]
\begin{center}
\begin{tabular}{cc}
\widgraph{.5\textwidth}{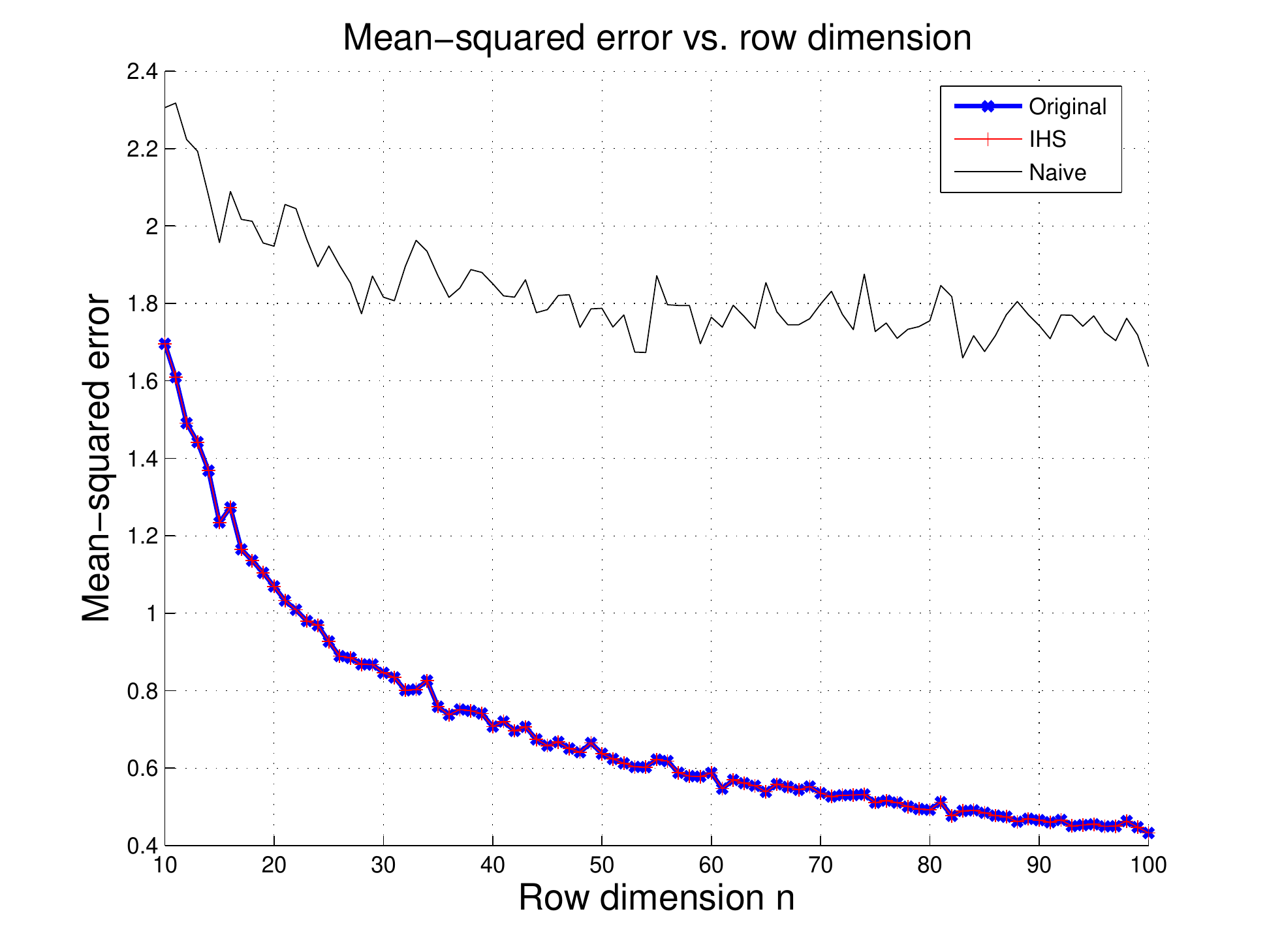} &
\widgraph{.51\textwidth}{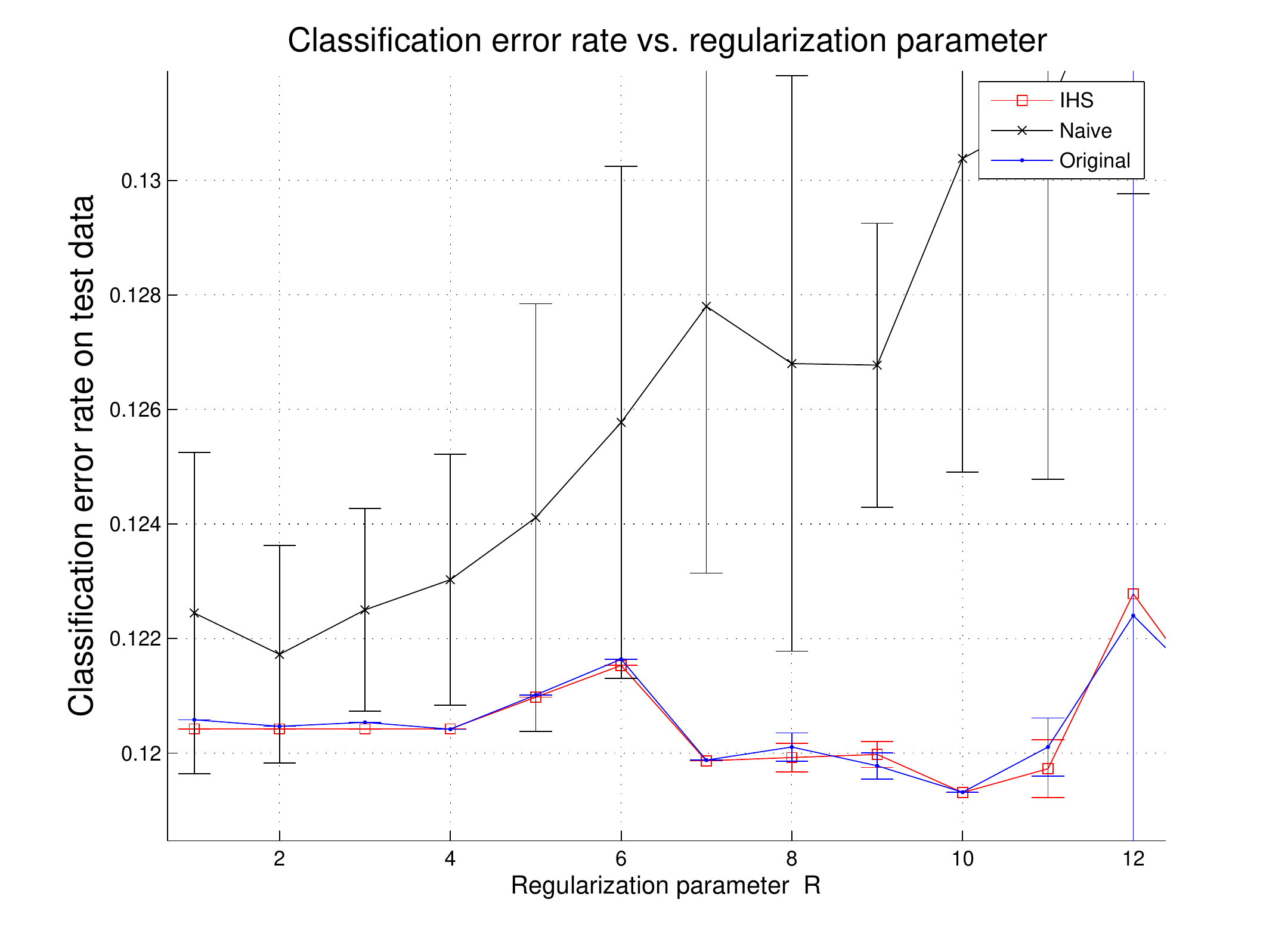} \\
(a) & (b)
\end{tabular}
\caption{Simulations of the IHS algorithm for nuclear-norm constrained
  problems. The blue curves correspond to the solution of the original
  (unsketched problem), wheras red curves correspond to the IHS method
  applied for $\Tcrit = 1 + \lceil \log(\numobs) \rceil$ rounds using
  a sketch size of $\numproj$. Black curves correpsond to the naive
  sketch applied using $M = N \numproj$ projections in total,
  corresponding to the same number used in all iterations of the IHS
  algorithm. (a) Mean-squared error versus the row dimension $\numobs
  \in [10,100]$ for recovering a $20 \times 20$ matrix of rank $\rdim
  2$, using a sketch dimension $\numproj = 60$.  Note how the accuracy
  of the IHS algorithm tracks the error of the unsketched solution
  over a wide range of $\numobs$, whereas the classical sketch has
  essentially constant error.  (b) Classification error rate versus
  regularization parameter $R \in\{1,\ldots,12\}$, with error bars
  corresponding to one standard deviation over the test set.
  Sketching algorithms were applied to the JAFFE face expression using
  a sketch dimension of $M = 100$ for the classical sketch, and $N =
  5$ iterations with $\numproj = 20$ sketches per iteration for the
  IHS algorithm. }
 \label{FigLowRank}
\end{center}
\end{figure}

More specifically, in order to verify the classification accuracy of
the classifier obtained by IHT algorithm, we solved the original
convex program, the classical sketch based on ROS sketches of
dimension $\numproj = 100$, and also the corresponding IHS algorithm
using ROS sketches of size $20$ in each of $5$ iterations.  In this
way, both the classical and IHS procedures use the same total number
of sketches, making for a fair comparison.  We repeated each of these
three procedures for all choices of the radius
$R\in\{1,2,3,\ldots,12\}$, and then applied the resulting classifiers
to classify images in the test dataset.  For each of the three
procedures, we calculated the classification error rate, defined as
the total number of mis-classified images divided by $\ntest \times
\dtask$.  Panel (b) of Figure~\ref{FigLowRank} plots the resulting
classification errors versus the regularization parameter. The error
bars correspond to one standard deviation calculated over the
randomness in generating sketching matrices.  The plots show that the
IHS algorithm yields classifiers with performance close to that given
by the original solution over a range of regularizer parameters, and
is superior to the classification sketch.  The error bars also show
that the IHS algorithm has less variability in its outputs than the
classical sketch.

%


\section{Discussion}
\label{SecDiscussion}

In this paper, we focused on the problem of solution approximation (as
opposed to cost approximation) for a broad class of constrained
least-squares problem.  We began by showing that the classical
sketching methods are sub-optimal, from an information-theoretic point
of view, for the purposes of solution approximation.  We then proposed
a novel iterative scheme, known as the iterative Hessian sketch, for
deriving $\varepsilon$-accurate solution approximations.  We proved a
general theorem on the properties of this algorithm, showing that the
sketch dimension per iteration need grow only proportionally to the
statistical dimension of the optimal solution, as measured by the
Gaussian width of the tangent cone at the optimum.  By taking
$\log(1/\varepsilon)$ iterations, the IHS algorithm is guaranteed to
return an $\varepsilon$-accurate solution approximation with
exponentially high probability.

In addition to these theoretical results, we also provided empirical
evaluations that reveal the sub-optimality of the classical sketch,
and show that the IHS algorithm produces near-optimal estimators.
Finally, we applied our methods to a problem of facial expression
using a multi-task learning model applied to the JAFFE face database.
We showed that IHS algorithm applied to a nuclear-norm constrained
program produces classifiers with considerably better classification
accuracy compared to the naive sketch.

There are many directions for further research, but we only list here
some of them. The idea behind iterative sketching can also be applied
to problems beyond minimizing a least-squares objective function
subject to convex constraints. An important class of such problems are
$\ell_p$-norm forms of regression, based on the convex program
\begin{align*}
\min_{x \in \real^\usedim} \|\Amat x-y\|_p^p \qquad \mbox{for some $p
  \in [1, \infty]$.}
\end{align*}
The case of $\ell_1$-regression ($p = 1$) is an important special
case, known as robust regression; it is especially effective for data
sets containing outliers~\cite{Huber73}. Recent
work~\cite{clarkson2013fast} has proposed to find faster solutions of
the $\ell_1$-regression problem using the classical sketch (i.e.,
based on $(\Sketch \Amat, \Sketch \yvec)$) but with sketching matrices
based on Cauchy random vectors.  Based on the results of the current
paper, our iterative technique might be useful in obtaining sharper
bounds for solution approximation in this setting as well.


\subsection*{Acknowledgements}

Both authors were partially supported by Office of Naval Research MURI
grant N00014-11-1-0688, and National Science Foundation Grants
CIF-31712-23800 and DMS-1107000. In addition, MP was supported by a
Microsoft Research Fellowship.


\appendix


\section{Proof of lower bounds}
\label{AppThmSubOptimal}

This appendix is devoted to the verification of
condition~\eqref{EqnHanaShouldSleep} for different model classes,
followed by the proof of Theorem~\ref{ThmSubOptimal}.


\subsection{Verification of condition~\eqref{EqnHanaShouldSleep}}
\label{AppLemHanaShouldSleep}

We verify the condition for three different types of sketches.

\paragraph{Gaussian sketches:}  First, let $\Sketch 
\in \real^{\numproj \times \numobs}$ be a random matrix with
i.i.d. Gaussian entries.  We use the singular value decomposition to
write \mbox{$\Sketch = U \Lambda V^T$} where both $U$ and $V$ are
orthonormal matrices of left and right singular vectors.  By rotation
invariance, the columns $\{v_i\}_{i=1}^\numproj$ are uniformly
distributed over the sphere $\Sphere{\numobs}$.  Consquently, we have
\begin{align}
\Exs_\Sketch \big[ \Sketch^T \big( \Sketch \Sketch^T)^{-1} \Sketch
  \big] & = \Exs \sum_{i=1}^\numproj v_i v_i^T =
\frac{\numproj}{\numobs} I_{\numobs},
\end{align}
showing that condition~\eqref{EqnHanaShouldSleep} holds with $\eta =
1$.


\paragraph{ROS sketches:}
In this case, we have $\Sketch = \sqrt{\numobs} P H D$, where $P \in
\real^{\numproj \times \numobs}$ is a random picking matrix (with each
row being a standard basis vector).  We then have $\Sketch \Sketch^T =
\numobs I_\numproj$ and also $\Exs_P[P^T P] = \frac{\numproj}{\numobs}
I_\numobs$, so that
\begin{align*}
\Exs_\Sketch[ \Sketch^T(\Sketch \Sketch^T)^{-1} \Sketch] & =
\Exs_{D,P} [D H^T P^T P H D ]= \Exs_{D} [D H^T
  (\frac{\numproj}{\numobs} I_\numobs) HD] = \frac{\numproj}{\numobs}
I_{\numobs},
\end{align*}
showing that the condition holds with $\eta = 1$.


\paragraph{Weighted row sampling:}

Finally, suppose that we sample $\numproj$ rows independently using a
distribution $\{p_j\}_{j=1}^\numobs$ on the rows of the data matrix
that is $\alpha$-balanced~\eqref{EqnAlphaBalanced}.  Letting
$\mathcal{R} \subseteq \{1, 2, \ldots, \numobs \}$ be the subset of
rows that are sampled, and let $N_j$ be the number of times each row
is sampled.  We then have 
\begin{align*}
\Exs \Big[\Sketch^T \big(\Sketch \Sketch^T)^{-1} \Sketch \Big] & =
\sum_{j \in \mathcal{R}} \Exs[e_j e_j^T] \; = \; D,
\end{align*}
where $D \in \real^{\numobs \times \numobs}$ is a diagonal matrix with
entries $D_{jj} = \mprob[j \in \mathcal{R}]$.  Since the trials are
independent, the $j^{th}$ row is sampled at least once in $\numproj$
trials with probability $q_j = 1 - (1-p_j)^m$, and hence
\begin{align*}
\Exs_\Sketch \big[ \Sketch^T \big( \Sketch \Sketch^T)^{-1} \Sketch
  \big] &= \diag \big ( \{ 1-(1-p_i)^m \}_{i=1}^\numproj \big) \;
\preceq \; \big(1-(1-p_\infty)^\numproj \big) I_\numobs \; \preceq
\numproj p_\infty,
\end{align*}
where $p_\infty = \max_{j \in [\numobs]} p_j$.  Consequently, as long
as the row weights are $\alpha$-balanced~\eqref{EqnAlphaBalanced} so
that $p_\infty \leq \frac{\alpha}{\numobs}$, we have
\begin{align*}
\opnorm{\Exs_\Sketch \big[ \Sketch^T \big( \Sketch \Sketch^T)^{-1}
    \Sketch\big]} & \leq \alpha \frac{\numproj}{\numobs}
\end{align*}
showing that condition~\eqref{EqnHanaShouldSleep} holds with $\eta =
\alpha$, as claimed.


\subsection{Proof of Theorem~\ref{ThmSubOptimal}}

Let $\{z^j\}_{j=1}^\PackNum$ be a $1/2$-packing of $\Constraint_0$ in
the semi-norm $\SEMI{\cdot}$ and for a fixed $\delta\in (0,1/4)$,
define $x^j = 4 \delta z^j$. We thus obtain a collection of vectors in
$\Constraint_0$ such that
\begin{align*}
2 \delta \le \frac{1}{\sqrt{n}}\|A(x^j - x^k) \|_2 \le 8\delta \qquad
\mbox{for all $j \neq k$.}
\end{align*}
Letting $J$ be a random index uniformly distributed over $\{1, \ldots,
\PackNum \}$, suppose that conditionally on $J = j$, we observe the
sketched observation vector $\Sketch \yvec = \Sketch \Amat x^j +
\Sketch \wvec$, as well as the sketched matrix $\Sketch \Amat$.
Conditioned on $J = j$, the random vector $\Sketch \yvec$ follows a
$\NORMAL(\Sketch \Amat x^j, \sigma^2 \Sketch \Sketch^T)$ distribution,
denoted by $\mprob_{x^j}$.  We let $\Vvar$ denote the resulting
mixture variable, with distribution $\frac{1}{\PackNum}
\sum_{j=1}^\PackNum \Prob_{x^j}$.

Consider the multiway testing problem of determining the index $J$
based on observing $\Vvar$.  With this set-up, a standard reduction in
statistical minimax (e.g.,~\cite{Birge87,Yu97}) implies that, for any
estimator $\xdagger$, the worst-case mean-squared error is lower bounded
as
\begin{align}
\label{EqnPeas}
\sup_{\xstar \in \Constraint} \Exs_{\Sketch, w} \SEMI{\xdagger -
  \xstar}^2 & \geq \delta^2 \inf_{\psi} \Prob[\psi(\Vvar)
  \neq J],
\end{align}
where the infimum ranges over all testing functions $\psi$.
Consequently, it suffices to show that the testing error is lower
bounded by $1/2$.

In order to do so, we first apply Fano's inequality~\cite{Cover}
conditionally on the sketching matrix $\Sketch$ to see that
\begin{align}
\label{EqnSpinach}
\mprob[\psi(\Vvar) \neq J] \; = \; \Exs_\Sketch \Big \{
\mprob[\psi(\Vvar) \neq J \mid \Sketch] \Big \} & \geq 1 - \frac{
  \Exs_\Sketch\big[I_\Sketch(\Vvar; J) \big] + \log 2} {\log
  \PackNum},
\end{align}
where $I_\Sketch(\Vvar; J)$ denotes the mutual information between
$\Vvar$ and $J$ with $\Sketch$ fixed.  Our next step is to upper bound
the expectation $\Exs_\Sketch[I(\Vvar; J)]$.

Letting $\kull{\Prob_{x^j}}{\mprob_{x^k}}$ denote the Kullback-Leibler
divergence between the distributions $\Prob_{x^j}$ and $\Prob_{x^k}$,
the convexity of Kullback-Leibler divergence implies that
\begin{align*}
I_\Sketch(\Vvar; J) & = \frac{1}{\PackNum} \sum_{j=1}^\PackNum
\kull{\Prob_{x^j}}{\frac{1}{\PackNum} \sum_{k=1}^\PackNum \Prob_{x^k}}
\; \leq \; \frac{1}{\PackNum^2} \sum_{j,k=1}^\PackNum
\kull{\mprob_{x^j}}{\mprob_{x^k}}.
\end{align*}
Computing the KL divergence for Gaussian vectors yields
\begin{align*}
I_\Sketch(\Vvar; J) & \leq \frac{1}{\PackNum^2} \sum_{j,k=1}^\PackNum
\frac{1}{2\sigma^2} (x^j-x^k)^T \Amat^T \Big[ \Sketch^T (\Sketch
  \Sketch^T)^{-1}\Sketch \Big] \Amat(x^j-x^k).
\end{align*}
Thus, using condition~\eqref{EqnHanaShouldSleep}, we have
\begin{align*}
\Exs_\Sketch[ I(\Vvar; J) ] & \leq \frac{1}{\PackNum^2
}\sum_{j,k=1}^\PackNum \frac{\numproj \; \eta}{2\, \numobs \sigma^2}
\|\Amat (x^j-x^k)\|_2^2 \; \leq \frac{32 \, \numproj \,
  \eta}{\sigma^2} \, \delta^2,
\end{align*}
where the final inequality uses the fact that $\SEMI{x^j- x^k} \leq 8
\delta$ for all pairs.

Combined with our previous bounds~\eqref{EqnPeas}
and~\eqref{EqnSpinach}, we find that
\begin{align*}
\sup_{\xstar \in \Constraint} \Exs \| \xhat - \xstar\|_2^2 & \geq
\delta^2 \Big \{ 1 - \frac{32 \frac{\numproj \, \eta \,
    \delta^2}{\sigma^2} + \log 2}{\log \PackNum} \Big \}.
\end{align*}
Setting $\delta = \frac{\sigma^2 \log(\PackNum/2)}{64 \, \eta \,
  \numproj }$ yields the lower bound~\eqref{EqnMinimax}.
%


\section{Proof of Proposition~\ref{PropHessSketch}}
\label{AppPropHessianSketch}

Since $\xhat$ and $\xls$ are optimal and feasible, respectively, for
the Hessian sketch program~\eqref{EqnHessianSketch}, we have
\begin{subequations}
\begin{align}
\inprod{\Amat^T \Sketch^T \big(\Sketch \Amat \xhat - y \big)}{\xls -
  \xhat} & \geq 0
\end{align}
Similarly, since $\xls$ and $\xhat$ are optimal and feasible, respectively,
for the original least squares program
\begin{align}
\label{EqnOptimalForOriginal}
\inprod{\Amat^T (\Amat \xls - \yvec)}{\xhat - \xls} & \geq 0.
\end{align}
\end{subequations}
Adding these two inequalities and performing some algebra yields the
basic inequality
\begin{align}
\label{EqnBasicInequality}
\frac{1}{\numproj} \|\Sketch \Amat \Delta\|_2^2 & \leq \Big| (\Amat
\xls)^T \big (I_\numobs - \frac{\Sketch^T \Sketch}{\numproj} \big)
\Amat \Delta \Big|.
\end{align}
Since $\Amat \xls$ is independent of the sketching matrix and $\Amat
\Delta \in \KCONELS$, we have
\begin{align*}
\frac{1}{\numproj} \|\Sketch \Amat \Delta\|_2^2 \geq \ZINF \, \|\Amat
\Delta\|_2^2, \qquad \mbox{and} \quad \Big| (\Amat \xls)^T \big
(I_\numobs - \Sketch^T \Sketch \big) \Amat \Delta \Big| \leq \ZSUP
\|\Amat \xls\|_2 \, \|\Amat \Delta\|_2,
\end{align*}
using the definitions~\eqref{EqnDefnZinf} and~\eqref{EqnDefnZsup} of
the random variables $\ZINF$ and $\ZSUP$ respectively.  Combining
the pieces yields the claim.


\section{Proof of Theorem~\ref{ThmOptIterative}}
\label{AppThmOptIterative}
It suffices to show that, for each iteration $t = 0, 1, 2, \ldots$, we
have
\begin{align}
\label{EqnInterBound}
\SEMI{ \xit{t+1} - \xls} & \leq \frac{\ZSUPIT{t+1}}{\ZINFIT{t+1}}
\SEMI{ \xit{t} - \xls}.
\end{align}
The claimed bounds~\eqref{EqnProdBound} and ~\eqref{EqnFinalBound}
then follow by applying the bound~\eqref{EqnInterBound} successively
to iterates $1$ through $\Tcrit$.

For simplicity in notation, we abbreviate $\SketchIt{t+1}$ to
$\Sketch$ and $\xit{t+1}$ to $\xhat$.  Define the error vector $\Delta
= \xhat - \xls$.  With some simple algebra, the optimization
problem~\eqref{EqnXitUpdate} that underlies the update $t+1$ can be
re-written as
\begin{align*}
\xhat & = \arg \min_{x \in \Constraint} \ENCMIN{ \frac{1}{2 \numproj}
  \|\Sketch \Amat x\|_2^2 - \inprod{\Amat^T \ytil}{x}},
\end{align*}
where $\ytil \defn \yvec + \Big [I - \frac{\Sketch^T
    \Sketch}{\numproj} \Big] \Amat \xit{t}$.  Since $\xhat$ and $\xls$
are optimal and feasible respectively, the usual first-order optimality
conditions imply that
\begin{align*}
\inprod{\Amat^T \frac{\Sketch^T \Sketch}{\numproj} \Amat x - \Amat^T
  \ytil}{\xls - \xhat} & \geq 0.
\end{align*}
As before, since $\xls$ is optimal for the original program, we have
\begin{align*}
\inprod{\Amat^T (\Amat \xls - \ytil + \Big [I - \frac{\Sketch^T
      \Sketch}{\numproj} \Big] \Amat \xit{t})}{\xhat - \xls} & \geq 0.
\end{align*}
Adding together these two inequalities and introducing the shorthand
$\Delta = \xhat - \xls$ yields
\begin{align}
\label{EqnRecursiveFinal}
\frac{1}{\numproj} \|\Sketch \Amat \Delta\|_2^2 & \leq \Big| (\Amat
(\xls - \xit{t})^T \big[ I - \frac{\Sketch^T \Sketch}{\numproj}
  \big] \Amat \Delta \Big|
\end{align}
Note that the vector $\Amat (\xls - \xit{t})$ is independent of the
randomness in the sketch matrix $\SketchIt{t+1}$.  Moreover, the
vector $\Amat \Delta$ belongs to the cone $\KCONE$, so that by the
definition of $\ZSUPIT{t+1}$, we have
\begin{subequations}
\begin{align}
\label{EqnRecursiveUpper}
\Big| (\Amat (\xls - \xit{t})^T \Big[ I - \frac{\Sketch^T
    \Sketch}{\numproj} \Big] \Amat \Delta \Big| & \leq \|\Amat (\xls -
\xit{t}) \|_2 \; \|\Amat \Delta\|_2 \; \ZSUPIT{t+1}.
\end{align}
Similarly, note the lower bound
\begin{align}
\label{EqnRecursiveLower}
\frac{1}{ \numproj} \|\Sketch \Amat \Delta\|_2^2 & \geq \|\Amat
\Delta\|_2^2 \; \ZINFIT{t+1}.
\end{align}
\end{subequations}
Combining the two bounds~\eqref{EqnRecursiveUpper}
and~\eqref{EqnRecursiveLower} with the earlier
bound~\eqref{EqnRecursiveFinal} yields the
claim~\eqref{EqnInterBound}.


\section{Maximum likelihood estimator and examples}
\label{AppPropMLE}

In this section, we a general upper bound on the error of the
constrained least-squares estimate.  We then use it (and other
results) to work through the calculations underlying
Examples~\ref{ExaOrdinaryLeastSquares} through~\ref{ExaMultiTask} from
Section~\ref{SecSubOpt}.


\subsection{Upper bound on MLE}

The accuracy of $\xls$ as an estimate of $\xstar$ depends on the
``size'' of the star-shaped set
\begin{align}
\label{EqnXstarCone}
\KCONE(\xstar) & = \big \{ v \in \real^\usedim \, \mid v =
\frac{t}{\sqrt{\numobs}} \, \Amat (x - \xstar) \quad \mbox{for some $t
  \in [0,1]$ and $x \in \Constraint$} \big \}.
\end{align}
When the vector $\xstar$ is clear from context, we use the shorthand
notation $\KCONE^*$ for this set.  By taking a union over all possible
$\xstar \in \Constraint_0$, we obtain the set $\KCONEUN \defn \bigcup
\limits_{\xstar \in \Constraint_0} \KCONE(\xstar)$, which plays an
important role in our bounds.  The complexity of these sets can be
measured of their \emph{localized Gaussian widths}.  For any radius
$\varepsilon > 0$ and set $\NewPlain \subseteq \real^\numobs$, the
Gaussian width of the set $\NewPlain \cap \Ball_2(\varepsilon)$ is
given by
\begin{subequations}
\begin{align}
\NonGaussComp_\varepsilon(\NewPlain) & \defn \Exs_g \Big[ \sup_{
    \substack{\theta \in \NewPlain \\ \|\theta\|_2 \leq \varepsilon}}
  |\inprod{w}{\theta}| \Big],
\end{align}
where $g \sim N(0, I_{\numobs \times \numobs})$ is a standard Gaussian
vector.  Whenever the set $\NewPlain$ is star-shaped, then it can be
shown that, for any $\sigma > 0$ and positive integer $\ell$, the
inequality
\begin{align}
\label{EqnCriticalValue}
\frac{\NonGaussComp_\varepsilon(\NewPlain)}{\varepsilon \,
  \sqrt{\ell}} & \leq \frac{\varepsilon}{\sigma}
\end{align}
\end{subequations}
has a smallest positive solution, which we denote by
$\DelCrit{\ell}(\Theta; \sigma)$.  We refer the reader to Bartlett et
al.~\cite{Bar05} for further discussion of such localized complexity
measures and their properties.

\noindent The following result bounds the mean-squared error
associated with the constrained least-squares estimate:
\bprops
\label{PropMLE}
For any set $\Constraint$ containing $\xstar$, the constrained
least-squares estimate~\eqref{EqnConstrainedLeastSquares} has
mean-squared error upper bounded as
\begin{align}
\label{EqnConstrainedMLEUpper}
\Exs_\wvec \big[\SEMI{\xls - \xstar}^2 \big] & \leq c_1 \big \{
\DelCrit{\numobs}^2 \big( \KCONESTAR \big) + \frac{\sigma^2}{\numobs}
\big \} \; \leq \; c_1 \big \{ \DelCrit{\numobs}^2 \big(\KCONEUN \big)
+ \frac{\sigma^2}{\numobs} \big \}.
\end{align}
\eprops

\noindent We provide the proof of this claim in
Section~\ref{AppProofPropMLE}.


\subsection{Detailed calculations for illustrative examples}
\label{AppExamples}

In this appendix, we collect together the details of calculations used
in our illustrative examples from Section~\ref{SecSubOpt}.  In all
cases, we make use tof the convenient shorthand $\Atil =
\Amat/\sqrt{\numobs}$.

\subsubsection{Unconstrained least squares:  
Example~\ref{ExaOrdinaryLeastSquares}}

By definition of the Gaussian width, we have
\begin{align*}
\NonGaussComp_\delta(\KCONESTAR) & = \Exs_g \big[ \sup_{\| \Atil \, (x
    - x^*) \|_2 \leq \delta} |\inprod{g}{\Atil (x - x^*)}| \big] \;
\leq \; \delta \, \sqrt{\usedim}
\end{align*}
since the vector $\Atil (x - \xstar)$ belongs to a subspace of
dimension $\rank(\Amat) = \usedim$.  The claimed upper
bound~\eqref{EqnUpperOLS} thus follows as a consequence of
Proposition~\ref{PropMLE}.


\subsubsection{Sparse vectors:  Example~\ref{ExaSparseLinear}}

The RIP property of order $8 \spindex$ implies
that
\begin{align*}
\frac{\|\Delta\|_2^2}{2} \; \stackrel{(i)}{\leq} \|\Atil \Delta\|_2^2
\; \stackrel{(ii)}{\leq} \; 2 \|\Delta\|_2^2 \qquad \mbox{for all
  vectors with $\|\Delta\|_0 \leq 8 \spindex$,}
\end{align*}
a fact which we use throughout the proof.
By definition of the Gaussian width, we have
\begin{align*}
\NonGaussComp_\delta(\KCONESTAR) & = \Exs_g \big[ \sup_{ \substack{
      \|x\|_1 \leq \|x^*\|_1 \\ \| \Atil (x - x^*) \|_2 \leq \delta}}
  |\inprod{g}{ \Atil ( x - x^*)}| \big].
\end{align*}
Since $x^* \in \Ball_0(\kdim)$, it can be shown (e.g., see the proof
of Corollary 3 in the paper~\cite{PilWai14a}) that for any vector
$\|x\|_1 \leq \|x^*\|_1$, we have $\|x - x^*\|_1 \leq 2 \sqrt{\kdim}
\|x - x^*\|_2$.  Thus, it suffices to bound the quantity
\begin{align*}
F(\delta; \kdim) & \defn \Exs_g \big[ \sup_{ \substack{ \|\Delta\|_1
      \leq 2 \sqrt{\kdim} \|\Delta\|_2 \\ \| \Atil \Delta \|_2 \leq
      \delta}} |\inprod{g}{ \Atil \Delta}| \big].
\end{align*}
By Lemma 11 in the paper~\cite{LohWai11}, we have $\Ball_1(\sqrt{s})
\cap \Ball_2(1) \subseteq 3 \clconv \Big \{ \Ball_0(s) \cap \Ball_2(1)
\Big \}$, where $\clconv$ denotes the closed convex hull.  Applying
this lemma with $s = 4 \kdim$, we have
\begin{align*}
F(\delta; \kdim) & \leq 3 \big[ \sup_{ \substack{ \|\Delta\|_0 \leq 4
      \kdim \\ \| \Atil \Delta \|_2 \leq \delta}} |\inprod{g}{ \Atil
    \Delta}| \big] \; \leq \; 3 \Exs \big[ \sup_{ \substack{
      \|\Delta\|_0 \leq 4 \kdim \\ \|\Delta \|_2 \leq 2\delta}}
  |\inprod{g}{ \Atil \Delta}| \big],
\end{align*}
using the lower RIP property (i).  By the upper RIP property, for
any pair of vectors $\Delta, \Delta'$ with $\ell_0$-norms at most
$4 \kdim$, we have 
\begin{align*}
\var \big(\inprod{g}{ \Atil \Delta} - \inprod{g}{\Atil \Delta'} \big)
& \leq 2 \|\Delta - \Delta'\|_2^2 \; = \; 2 \var
\big(\inprod{g}{\Delta - \Delta'} \big)
\end{align*}
Consequently, by the Sudakov-Fernique comparison~\cite{LedTal91}, we have
\begin{align*}
\Exs \big[ \sup_{ \substack{ \|\Delta\|_0 \leq 4 \kdim \\ \|\Delta
      \|_2 \leq 2\delta}} |\inprod{g}{ \Atil \Delta}| \big] & \leq 2
\Exs \big[ \sup_{ \substack{ \|\Delta\|_0 \leq 4 \kdim \\ \|\Delta
      \|_2 \leq 2\delta}} |\inprod{g}{\Delta}| \big] \; \leq c \;
\delta \sqrt{\kdim \log \big(\frac{e \usedim}{\kdim} \big)},
\end{align*}
where the final inequality standard results on Gaussian
widths~\cite{Gor07}.  All together, we conclude that
\begin{align*}
\DelCrit{\numobs}^2(\KCONESTAR; \sigma) & \leq c_1 \sigma^2
\frac{\kdim \log \big( \frac{e \usedim}{\kdim} \big)}{\numobs}.
\end{align*}
Combined with Proposition~\ref{PropMLE}, the claimed upper
bound~\eqref{EqnUpperSparse} follows.

 In the other direction, a straightforward argument
 (e.g.,~\cite{RasWaiYu09}) shows that there is a universal constant $c
 > 0$ such that $\log \PackNum_{1/2} \geq c \, \kdim \log \big(\frac{e
   \usedim}{\kdim} \big)$, so that the stated lower bound follows from
 Theorem~\ref{ThmSubOptimal}.


\subsubsection{Low rank matrices:  Example~\ref{ExaMultiTask}:}

By definition of the Gaussian width, we have width, we have
\begin{align*}
\NonGaussComp_\delta(\KCONESTAR) & = \Exs_g \Biggr[ \sup_{
    \substack{\fronorm{\Atil \, (X - \Xstar)} \leq \delta
      \\ \nucnorm{X} \leq \nucnorm{\Xstar}}} |\tracer{\Atil^T G}{(X -
    \Xstar)}| \Biggr],
\end{align*}
where $G \in \real^{\numobs \times \usedimb}$ is a Gaussian random
matrix, and $\tracer{C}{D}$ denotes the trace inner product between
matrices $C$ and $D$.  Since $\Xstar$ has rank at most $\rdim$, it can
be shown that $\nucnorm{X - \Xstar} \leq 2 \sqrt{\rdim} \fronorm{X -
  \Xstar}$; for instance, see Lemma 1 in the paper~\cite{NegWai09}.
Recalling that $\gammin(\Atil)$ denotes the minimum singular value, we
have
\begin{align*}
\fronorm{X - \Xstar} \leq \frac{1}{\gammin(\Atil)} \, \fronorm{\Atil
  (X - \Xstar)} \leq \frac{\delta}{\gammin(\Atil)}.
\end{align*}
Thus, by duality between the nuclear and operator norms, we have
\begin{align*}
 \Exs_g \Biggr[ \sup_{ \substack{\fronorm{\Atil \, (X - \Xstar)} \leq
       \delta \\ \nucnorm{X} \leq \nucnorm{\Xstar}}} |\tracer{G}{\Atil
     (X - \Xstar)}| \Biggr] & \leq \frac{2 \, \sqrt{\rdim} \,
   \delta}{\gammin(\Amat)} \, \Exs \big[ \opnorm{\Atil^T G}].
\end{align*}
Now consider the matrix $\Amat^T G \in \real^{\usedima \times
  \usedimb}$.  For any fixed pair of vectors \mbox{$(u, v) \in
  \Sphere{\usedima} \times \Sphere{\usedimb}$,} the random variable $Z
= u^T \Atil^T G v$ is zero-mean Gaussian with variance at most
$\gammax^2(\Atil)$.  Consequently, by a standard covering argument in
random matrix theory~\cite{Ver11}, we have $\Exs \big[ \opnorm{\Atil^T
    G}] \precsim \gammax(\Atil) \big( \sqrt{\usedima + \usedimb}
\big)$.  Putting together the pieces, we conclude that
\begin{align*}
\DelCrit{\numobs}^2 & \preceq \sigma^2 \;
\frac{\gammax^2(\Amat)}{\gammin^2(\Amat)} \, \rdim \, (\usedima +
\usedimb),
\end{align*}
so that the upper bound~\eqref{EqnUpperLowRank} follows from
Proposition~\ref{PropMLE}.

\subsection{Proof of Proposition~\ref{PropMLE}}
\label{AppProofPropMLE}

Throughout this proof, we adopt the shorthand $\deln =
\DelCrit{\numobs}(\KCONESTAR)$.  Our strategy is to prove the
following more general claim: for any $t \geq \deln$, we have
\begin{align}
\label{EqnMSETail}
\mprob_{\Sketch, \wvec} \big[ \SEMI{\xls - \xstar}^2 \geq 16 t \deln
  \big] & \leq c_1 \CEXP{- c_2 \frac{\numobs t \deln}{\sigma^2}}.
\end{align}
A simple integration argument applied to this tail bound implies the
claimed bound~\eqref{EqnConstrainedMLEUpper} on the expected
mean-squared error.

Since $\xstar$ and $\xls$ are feasible and optimal, respectively, for
the optimization problem~\eqref{EqnConstrainedLeastSquares}, we have
the basic inequality
\begin{align*}
\frac{1}{2 \numobs} \|y - A \xls\|_2^2 & \leq \frac{1}{2 \numobs} \|y
- \Amat \xstar\|_2 \; = \; \frac{1}{2 \numobs} \|w\|_2^2.
\end{align*}
Introducing the shorthand $\Delta = \xls - \xstar$ and re-arranging
terms yields
\begin{align}
\label{EqnMLEBasic}
\frac{1}{2} \SEMI{\Delta}^2 \; = \; \frac{1}{2\numobs} \|A
\Delta\|_2^2 & \leq \frac{\sigma}{\numobs} \big| \sum_{i=1}^\numobs
\inprod{g}{A \Delta} \big|,
\end{align}
where $g \sim N(0, I_\numobs)$ is a standard normal vector.

For a given $u \geq \deln$, define the ``bad'' event
\begin{align*} 
\AuxEvent(u) & \defn \big \{ \exists \quad z \in \Constraint - \xstar
\quad \mbox{with $\SEMI{z} \geq u$, and $|\frac{\sigma}{\numobs}
  \sum_{i=1}^\numobs g_i (\Amat z)_i |\geq 2 u \, \SEMI{z} $} \big \}
\end{align*} 
The following lemma controls the probability of this event:
\blems
\label{LemCore}
For all $u \geq \deln$, we have $\mprob[\AuxEvent(u)] \leq
\CEXP{-\frac{ \numobs u^2}{2 \sigma^2}}$.
\elems

Returning to prove this lemma momentarily, let us prove the
bound~\eqref{EqnMSETail}.  For any $t \geq \deln$, we can apply
Lemma~\ref{LemCore} with $u = \sqrt{t \deln}$ to find that
\begin{align*}
\mprob[\AuxEvent^c(\sqrt{t \deln})] & \geq 1 - \CEXP{-\frac{
    \numobs t \deln}{2 \sigma^2}}.
\end{align*}

If $\SEMI{\Delta} < \sqrt{t \, \deln}$, then the claim is immediate.
Otherwise, we have $\SEMI{\Delta} \geq \sqrt{ t \, \deln}$.  Since
$\Delta \in \Constraint - \xstar$, we may condition on
$\AuxEvent^c(\sqrt{t \deln})$ so as to obtain the bound
\begin{align*}
\big |\frac{\sigma}{\numobs} \sum_{i=1}^\numobs g_i (\Amat \Delta)_i
\big| & \leq 2 \, \SEMI{\Delta} \: \sqrt{t \deln}.
\end{align*}
Combined with the basic inequality~\eqref{EqnMLEBasic}, we see that
\begin{align*}
\frac{1}{2} \SEMI{\Delta}^2 \leq 2 \, \SEMI{\Delta} \: \sqrt{t \deln},
\qquad \mbox{or equivalently $\SEMI{\Delta}^2 \leq 16 t \deln$,}
\end{align*}
a bound that holds with probability greater than $1 -
\CEXP{-\frac{\numobs t \deln}{2 \sigma^2}}$ as claimed.\\

It remains to prove Lemma~\ref{LemCore}.  Our proof involves the
auxiliary random variable
\begin{align*}
\MyMax(u) \, \defn \, \sup_{\substack{z \in \starset(\Constraint -
    \xstar) \\ \SEMI{z} \leq u}} |\frac{\sigma}{\numobs}
\sum_{i=1}^\numobs g_i \: (\Amat z)_i|,
\end{align*}

\paragraph{Inclusion of events:}  We first claim that
$\AuxEvent(u) \subseteq \{ \MyMax(u) \geq 2 u^2 \}$. Indeed, if
$\AuxEvent(u)$ occurs, then there exists some $z \in \Constraint -
\xstar$ with $\SEMI{z} \geq u$ and
\begin{align}
\label{EqnKeyViolate}
|\frac{\sigma}{\numobs} \sum_{i=1}^\numobs g_i \: (\Amat z)_i| & \geq
2 u \; \SEMI{z}.
\end{align}
Define the rescaled vector $\ztil = \frac{u}{\SEMI{z}} z$.  Since $z
\in \Constraint - \xstar$ and $\frac{u}{\SEMI{z}} \leq 1$, the vector
$\ztil \in \starset(\Constraint - \xstar)$.  Moreover, by
construction, we have $\SEMI{\ztil} = u$.  When the
inequality~\eqref{EqnKeyViolate} holds, the vector $\ztil$ thus
satisfies $|\frac{\sigma}{\numobs} \sum_{i=1}^\numobs g_i \: (\Amat
\ztil)_i| \geq 2 u^2$, which certifies that $\MyMax(u) \geq 2 u^2$, as
claimed.

\paragraph{Controlling the tail probability:}  The final step
is to control the probabability of the event $\{\MyMax(u) \geq 2
u^2\}$.  Viewed as a function of the standard Gaussian vector $(g_1,
\ldots, g_\numobs)$, it is easy to see that $\MyMax(u)$ is Lipschitz
with constant $L = \frac{\sigma u}{\sqrt{\numobs}}$.  Consequently, 
by concentration of measure for Lipschitz Gaussian functions, we have
\begin{align}
\label{EqnMyVarConcentrate}
\mprob \big[ \MyMax(u) \geq \Exs[\MyMax(u)] + u^2 \big] & \leq \CEXP{-
  \frac{\numobs u^2}{2 \sigma^2}}.
\end{align}
In order to complete the proof, it suffices to show that
$\Exs[\MyMax(u)] \leq u^2$.  By definition, we have $\Exs[\MyMax(u)] =
  \frac{\sigma}{\sqrt{\numobs}} \NonGaussComp_u(\KCONESTAR)$.  Since
  $\KCONESTAR$ is a star-shaped set, the function $v \mapsto
  \NonGaussComp_v(\KCONESTAR)/v$ is non-increasing~\cite{Bar05}.
Since $u \geq \deln$, we have
\begin{align*}
\sigma \frac{\NonGaussComp_u(\KCONESTAR)}{u} & \leq \sigma
\frac{\NonGaussComp_{\deln}(\KCONESTAR)}{\deln} \; \leq \; \deln.
\end{align*}
where the final step follows from the definition of $\deln$.  Putting
together the pieces, we conclude that $\Exs[\MyMax(u)] \leq \deln u
\leq u^2$ as claimed.


\bibliographystyle{plain}

\bibliography{P_OptimalIterative}
\end{document}